\documentstyle[12pt,amsmath,amsfonts,amssymb,leqno]{article}
\newtheorem{theorem}{Theorem}[section]
\newtheorem{definition}[theorem]{Definition}
\newtheorem{corollary}[theorem]{Corollary}
\newtheorem{proposition}[theorem]{Proposition}
\newtheorem{remark}[theorem]{Remark}
\newtheorem{lemma}[theorem]{Lemma}
\newtheorem{claim}[theorem]{Claim}
\newtheorem{example}[theorem]{Example}
\newtheorem{problem}[theorem]{Problem}
\newcommand {\Mcc}     {{\mathcal M}}

\newcommand {\Pc}      {{\cal P}}
\newcommand {\Kc}      {{\cal K}}

\newcommand {\Cc}      {{\cal C}}
\newcommand {\Hc}      {{\cal H}}
\newcommand {\Tc}      {{\cal T}}
\newcommand {\Gc}      {{\cal G}}

\newcommand {\R}       {{\bf R}}
\newcommand {\RN}      {\R^n}
\newcommand {\KS}      {\Kc(S)}

\newcommand {\HP}      {{\bf H}^{n+1}}
\newcommand {\KRN}     {\Kc}

\newcommand {\RNP}     {{\bf R}^{n+1}_+}
\newcommand {\TOM}     {\Tc_\omega}

\newcommand {\PRMP}    {\PM\times\RNP}
\newcommand {\PMKS}    {\PM\times\KS}
\newcommand {\PMK}     {\PM\times\Kc}

\newcommand {\rp}      {\varrho}
\newcommand {\rph}     {\rho_{H}}

\newcommand {\PK}      {\Pc_{k}}
\newcommand {\PM}      {\Pc_{m-1}}
\newcommand {\PL}      {\Pc_{L}}

\newcommand {\PKR}     {\Pc_{L}\times\KRN}

\newcommand {\MRR}     {(\Mcc,\rho)}

\newcommand {\COR}     {C^{k,\omega}(\R)}

\newcommand {\CKO}     {C^{k,\omega }(\R^{n})}
\newcommand {\COM}     {C^{1,\omega }(\R^{n})}

\newcommand {\CK}      {C^{k}(\RN)}

\newcommand {\CKLM}    {C^k\Lambda^m_{\omega }(\R^{n})}
\newcommand {\LMO}     {\Lambda^m_{\omega }(\R^{n})}
\newcommand {\LMT}     {\Lambda^2_{\omega }(\R^{n})}
\newcommand {\LZ}      {{\bf Z}(\R^{n})}
\newcommand {\ZM}      {{\bf Z}_m(\R^{n})}
\newcommand {\card}    {\operatorname{card}}
\newcommand {\Lip}     {\operatorname{Lip}}

\newcommand {\BLO}     {\operatorname{{\bf LO}}}

\newcommand {\dd}      {\operatorname{d}}

\newcommand {\dm}      {\dd_\omega}
\newcommand {\dom}     {\delta_{\omega}}
\newcommand {\tP}      {\widetilde{P}}

\newcommand {\tQ}      {\widetilde{Q}}

\newcommand {\tr}      {\tilde{r}}
\newcommand {\tx}      {\tilde{x}}

\newcommand {\br}      {\bar{r}}

\newcommand {\bQ}      {K}

\newcommand {\fa}      {\varphi_{\alpha}}
\newcommand {\fin}     {\varphi_{\alpha}^{-1}}

\newcommand {\domg}    {\delta_\omega}
\newcommand {\rom}     {\rho_\omega}
\newcommand {\nn}      {\nonumber}
\newcommand {\reff}[1] {\ref{#1}}
\newcommand {\rf}[1]   {(\ref{#1})}

\newcommand {\bx}      {$\hspace*{15mm}\Box$}

\newcommand {\SECT}[2]
{\section*{\centerline{\normalsize {\bf #1}}}
\setcounter{section}{#2}
\setcounter{theorem}{0}\setcounter{equation}{0}}
\newcommand{\lbl}[1]   {\label{#1}}
\newcommand{\be}       {\begin{eqnarray}}
\newcommand{\bel}[1]   {\begin{eqnarray}\label{#1}}
\newcommand{\ee}       {\end{eqnarray}}
\begin{document}
\medskip
\centerline{\large{\bf The Whitney extension problem for
Zygmund spaces}} \vspace*{5mm}
\centerline{\large{\bf and Lipschitz selections in
hyperbolic jet-spaces}}
\vspace*{8mm}
\centerline{By {\it Pavel Shvartsman}} \vspace*{5 mm}
\centerline {\it Department of Mathematics, Technion -
Israel Institute of Technology,}\vspace*{2 mm}
\centerline{\it 32000 Haifa, Israel}\vspace*{2 mm}
\centerline{\it e-mail: pshv@tx.technion.ac.il}
\vspace*{10 mm}
\renewcommand{\thefootnote}{ }
\footnotetext[1]{{\it\hspace{-6mm}Math Subject
Classification} 46E35, 52A35, 54C60, 54C65 }
\begin{abstract} {\small \par We study a variant of the
Whitney extension problem \cite{W1,W2} for the space
$\CKLM$ of functions whose partial derivatives of order $k$
satisfy the generalized Zygmund condition. We identify
$\CKLM$ with a space of {\it Lipschitz} mappings from a
metric space  $(\RNP,\rom)$ supplied with a hyperbolic
metric $\rom$ into a metric space
$(\Pc_{k+m-1}\times\RNP,\dm)$ of polynomial fields on
$\RNP$ equipped with a hyperbolic-type metric $\dm$. This
identification allows us to reformulate the Whitney problem
for $\CKLM$ as a Lipschitz selection problem for set-valued
mappings from $(\RNP,\rom)$ into a certain family of
subsets of $\Pc_{k+m-1}\times\RNP$.}
\end{abstract} \vspace*{15mm}
\renewcommand{\thefootnote}{\arabic{footnote}}
\setcounter{footnote}{0}
\SECT{1. Introduction}{1}
\indent
\par Let $m$ be a non-negative integer. We let $\Omega_m$
denote the class of non-decreasing continuous functions
$\omega:\R_+\to\R_+$ such that $\omega(0)=0$ and the
function $\omega(t)/t^m$ is non-increasing. Given
non-negative integers $k$ and $m$ and $\omega\in\Omega_m$
we define the space $\CKLM$ as follows: a function $f\in
\CK$ belongs to $\CKLM$ if there exists a constant
$\lambda>0$ such that for every multi-index $\alpha$,
$|\alpha|=k$, and every $x,h\in\RN$
$$
|\Delta^m_h(D^{\alpha}f)(x)|\le \lambda\omega(\|h\|).
$$
Here as usual $\Delta^m_hf$ denotes the $m$-th difference
of a function $f$ of step $h$, i.e., the quantity
$$
\Delta^m_hf(x):= \sum^m_{i=0}(-1)^{m-i}{m\choose
i}f(x+ih).
$$
\par $\CKLM$ is normed by
\bel{14.N} \|f\|_{\CKLM}:= \sum_{|\alpha|\le
k}\sup_{x\in\RN}|D^{\alpha}f(x)|+ \sum_{|\alpha|=
k}\sup_{x,h\in\RN}\frac{|\Delta^m_h(D^{\alpha}f)(x)|}
{\omega(\|h\|)}. \ee
\par In particular, for $m=1$ and $\omega\in\Omega_1$ the
space $C^k\Lambda^1_\omega(\RN)$ coincides with the space
$\CKO$ consisting of all functions $f\in \CK$ whose partial
derivatives of order $k$ satisfy the Lipschitz condition
(with respect to $\omega$):
$$
|D^{\alpha}f(x)-D^{\alpha}f(y)|\le\omega(\|x-y\|),
~~x,y\in\RN.
$$
In turn, the space $\LMO:=
C^0\Lambda^m_\omega(\RN),~\omega\in\Omega_m,$ coincides
with the generalized Zygmund space of bounded functions $f$
on $\RN$ whose modulus of smoothness of order $m$,
$\omega_m(\cdot;f)$, satisfies the inequality
$$
\omega_m(t;f)\le \lambda \omega(t),~~t\ge 0.
$$
In particular,  the space $\LMT$ with $\omega(t)=t$ is the
classical Zygmund space $\LZ$ of bounded functions
satisfying the Zygmund condition: there is $\lambda>0$ such
that for all $x,y\in\RN$
$$
|f\left(x\right)-2f\left(\tfrac{x+y}{2}\right)
+f(y)|\le\lambda \|x-y\|.
$$
(See, e.g. Stein \cite{St}.)
\par Throughout the paper we let $S$ denote an arbitrary
closed subset of $\RN$.
\par In this paper we study the following extension
problem.
\begin{problem}\lbl{WP-CKLM} Given non-negative integers
$k$ and $m$, a function $\omega\in \Omega_m$, and an
arbitrary function $f:S\to\R$, what is a necessary and
sufficient condition for $f$ to be the restriction to $S$
of a function $F\in\CKLM$?
\end{problem}
\par This is a variant of a classical problem which is
known in the literature as the Whitney Extension Problem
\cite {W1,W2}. It has attracted a lot of attention in
recent years. We refer the reader to \cite{BS1}-\cite{BS4},
\cite{F1}-\cite{F7}, \cite{BMP1,BMP2} and \cite{Z1,Z2} and
references therein for numerous results in this direction,
and for a variety of techniques for obtaining them.
\par  This note is devoted to the {\it phenomenon of
``finiteness"} in the Whitney problem for the spaces
$\CKLM$. It turns out that, in many cases, Whitney-type
problems for different spaces of smooth functions can be
reduced to the same kinds of problems, but for finite sets
with {\it prescribed numbers of points.}
\par For the space $\LMT$, this phenomenon has been studied
in the author's papers \cite{S0,Sf,S1}. It was shown that
{\it a function $f$ defined on $S$ can be extended to a
function $F\in \LMT$ with $\|F\|_{\LMT}\le
\gamma=\gamma(n)$ provided its restriction $f|_{S'}$ to
every subset $S'\subset S$ consisting of at most
$N(n)=3\cdot 2^{n-1}$ points can be extended to a function
$F_{S'}\in \LMT$ with $\|F_{S'}\|_{\LMT}\le 1$.} (Moreover,
the value $3\cdot 2^{n-1}$ is sharp \cite{S1}.)
\par This result is an example of {\it ``the finiteness
property"} of the space $\LMT$. We call the number $N$
appearing in formulations of finiteness properties {\it
``the finiteness number".}
\par In his pioneering work \cite{W2}, Whitney
characterized the restriction of the space $C^k(\R), k\ge
1,$ to an arbitrary subset $S\subset\R$ in terms of divided
differences of functions. An application of Whitney's
method to the space $\COR$ implies the finiteness property
for this space with the finiteness number $N=k+2$.
\par The restriction of the space
$C^k\Lambda^m_{\omega}(\R)$ to an arbitrary subset
$S\subset\R$ has been characterized by Jonsson \cite{J}
($m$ is arbitrary, $k=0$, $\omega(t)=t^{m-1}$), Shevchuk
\cite{She,She1} ($m,\omega$ are arbitrary, $k=0$), Galan
\cite{Gal} (the general case). These results imply the
finiteness property for $C^k\Lambda^m_{\omega}(\R)$ with
the finiteness number $N=m+k+1$.
\par For the space $\COM$ the finiteness property (with the
same finiteness number $N(n)=3\cdot 2^{n-1}$) has been
proved in \cite{BS4}, see also \cite{BS1}.
\par An impressive breakthrough in the solution of the
Whitney problem for $C^{k,\omega}$-spaces has recently been
made by Fefferman \cite{F1}-\cite{F7}. In particular, one
of his remarkable results states that the space $\CKO$
possesses the finiteness property {\it for all $k,n>1$,}
see \cite{F1,F3}. (An upper bound for the finiteness number
$N(k,n)$  is $N(k,n)\le 2^{\dim\PK}$, see Bierstone, Milman
\cite{BM}, and Shvartsman \cite{S5}. Here $\PK$ stands for
the space of polynomials of degree at most $k$ defined on
$\RN$. Recall that $\dim\PK={n+k\choose k}$.)
\par Thus for  $m>2$ or $m=2$ and $k>0$ the following
problem is open.
\begin{problem}\lbl{FP-CKLM} Whether the space $\CKLM$
possesses the finiteness property?
\end{problem}
\par In this paper we develop an approach to Problem
\reff{WP-CKLM} which allows us to reformulate this problem
as a purely geometric question about the existence of
Lipschitz selections of set-valued mappings defined on
metric spaces with certain hyperbolic structure.
\SECT{2. Extensions of Zygmund functions and Lipschitz
selections}{2}
\indent \par We will demonstrate this approach for the case
of the Zygmund space $\ZM:=C^0\Lambda^m_\omega(\RN)$ where
$\omega(t)=t^{m-1}$. Thus $\ZM$ is defined by the
finiteness of the norm
\bel{NR-ZM}
\|f\|_{\ZM}:= \sup_{x\in\RN}|f(x)|+
\sup_{x,h\in\RN}\frac{|\Delta^m_hf(x)|} {\|h\|^{m-1}}.
\ee
\par The crucial ingredient of our approach is an
isomorphism between the space $\ZM|_S$ and a space of
Lipschitz mappings from the set $S\times \R_+$ into the
product $\PRMP$ equipped with certain metrics.
\par This isomorphism is motivated by a description of the
restrictions of $C^k\Lambda^m_{\omega}$-functions in terms
of local approximations which we present in Section 3. Let
us formulate this result for the space  $\ZM$.
\par We will assume that all cubes in this paper are closed
and have sides which are parallel to the coordinate axes.
It will be convenient for us to measure distances in $\RN$
in the uniform norm
$$ \|x\|:=\max\{|x_i|:~i=1,...,n\}, \ \ \
x=(x_1,...,x_n)\in\RN. $$
Thus every cube
$$ Q=Q(x,r):=\{y\in\RN:\|y-x\|\le r\} $$
is a ``ball" in the metric space $(\RN,\|\cdot\|)$ of
``radius" $r$ centered at $x$. We let $x_Q:=x$ denote
center of $Q$ and $r_Q:=r$ its ``radius". Given a constant
$\lambda> 0$, we let $\lambda Q$ denote the cube
$Q(x,\lambda r)$.
\par We let
$$
\KS:=\{Q(x,r):~x\in S,~r>0\}
$$
denote the family of all cubes centered in $S$. By $\KRN$
we denote the family of {\it all cubes in $\RN$}; thus
$\KRN=\Kc(\RN)$.
\begin{theorem}\lbl{ZM-LOCAP}
(A) Let $f:S\to\R$. Suppose that there exists a function
$F\in\ZM$ such that $F|_S=f$. Then there exists a constant
$0<\lambda\le C\|F\|_{\CKLM}$ and a family of polynomials
$\{P_Q\in\PM:~Q\in\KS\}$ such that:
\par (1).  For every $x\in S$ and every cube $Q=Q(x,r),
r>0,$ we have $ P_Q(x)=f(x); $
\par (2). For every $Q\in\KS$ with $r_Q\le 1$ and every
$\beta, |\beta|\le m-1$,
\bel{C2} |D^{\beta} P_Q(x_Q)|\le \lambda\,r_Q^{-|\beta|}~;
\ee
\par (3). For every $Q_1=Q(x_1,r_1),Q_2=Q(x_2,r_2)\in\KS$,
and every $\alpha,|\alpha|\le m-1,$ we have
\be |D^{\alpha} P_{Q_1}(x_1)-D^{\alpha} P_{Q_2}(x_1)|&\le&
\lambda(\max\{r_1,r_2\}+\|x_1-x_2\|)^{m-1-|\alpha|}
\label{SE3}\\
&\cdot& \ln\left(1+
\frac{\max\{r_1,r_2\}+\|x_1-x_2\|}{\min\{r_1,r_2\}}\right)
.\nn \ee
\par (B) Conversely, suppose that there exists a constant
$\lambda>0$ and a family of polynomials
$\{P_Q\in\PM:~Q\in\KS\}$ such that conditions (2) and (3)
are satisfied.
Then for every $x\in S$ there exists the limit
\bel{14-ZM} f(x)=\lim_{x_Q=x,\,r_Q\to 0}P_Q(x). \ee
Moreover, there exists $F\in\ZM$ with $\|F\|_{\ZM}\le
C\lambda$ such that $F|_S=f$.
\par Here $C$ is a constant depending only on $m$ and $n$.
\end{theorem}
\par This result is a particular case of Theorems
\reff{T14.2} and \reff{T14.7} proven in Section 3. It can
be considered as a certain version of the classical Whitney
 extension theorem \cite{W1} for the Zygmund space $\ZM$
where the Taylor polynomials are replaced by corresponding
approximation polynomials $P_Q, Q\in\KS$.
\par Now our aim is to transform inequalities \rf{SE3} into
a certain Lipschitz condition for the mapping $\KS\ni
Q\to(P_Q,Q)\in\PMKS$. This is the crucial point of the
approach.
\par  We equip the family $\Kc$ (of all cubes in $\RN$)
with the distance $\rho:\Kc\times\Kc\to\R_+$ defined by the
following formula: if $Q_1=Q(x_1,r_1), Q_2=Q(x_2,r_2)\in
\Kc,$ and $Q_1\ne Q_2$, then
\bel{MHP}
\rho(Q_1,Q_2):=\ln\left(1+\frac{\max(r_1,r_2)+
\|x_1-x_2\|}{\min(r_1,r_2)}\right),
\ee
and $\rho(Q_1,Q_2):=0$ whenever $Q_1=Q_2.$ We prove that
$\rho$ is a metric on $\Kc$; moreover, the metric space
$(\Kc,\rho)$ can be identified (up to a constant weight)
with the classical Poin\^{c}are upper half-space model of
the hyperbolic space $H_{n+1}$, see Remark \reff{REMH}.
\par Now, for every $\alpha, |\alpha|\le m-1$, we will
rewrite every inequality in \rf{SE3} in such a way that its
right-hand side will be precisely equal  to
$\rho(Q_1,Q_2)$. By \rf{SE3}, we have
\be \frac{1}{\lambda}\,\,\frac{|D^{\alpha}
P_{Q_1}(x_1)-D^{\alpha}
P_{Q_2}(x_1)|}{\min\{r_1,r_2\}^{m-1-|\alpha|}}&\le&
\left(\frac{\max\{r_1,r_2\}+\|x_1-x_2\|}
{\min\{r_1,r_2\}}\right)^{m-1-|\alpha|}
\nn\\
&\cdot& \ln\left(1+
\frac{\max\{r_1,r_2\}+\|x_1-x_2\|}{\min\{r_1,r_2\}}\right)
\nn\\
&=&\left(e^{\rho(Q_1,Q_2)}-1\right)^{m-1-|\alpha|}\,
\rho(Q_1,Q_2). \nn \ee
\par Put
$$
\psi_\alpha(t):=t(e^t-1)^{m-1-|\alpha|},
$$
and by $\psi_\alpha^{-1}$ denote the inverse to
$\psi_\alpha$. Then, by the latter inequality,
$$
\frac{1}{\lambda}\,\,\frac{|D^{\alpha}
P_{Q_1}(x_1)-D^{\alpha}
P_{Q_2}(x_1)|}{\min\{r_1,r_2\}^{m-1-|\alpha|}}\le
\psi_\alpha(\rho(Q_1,Q_2)),
$$
so that
$$
\psi_\alpha^{-1}\left(\frac{1}{\lambda}
\,\,\frac{|D^{\alpha}
P_{Q_1}(x_1)-D^{\alpha}
P_{Q_2}(x_1)|}{\min\{r_1,r_2\}^{m-1-|\alpha|}}\right)\le
\rho(Q_1,Q_2).
$$
Of course, the same inequality holds for $x_2$ instead of
$x_1$. Taking the maximum over $x_1,x_2,$ and over all
$\alpha$ with $|\alpha|\le m-1$, we obtain
\bel{SEL}
I\left(\tfrac{1}{\lambda}\,\,P_{Q_1},
\tfrac{1}{\lambda}\,\,P_{Q_2}\right)\le \rho(Q_1,Q_2)
\ee
where
\bel{DFI} I(P_{Q_1},P_{Q_2}):=\max_{i=1,2,|\alpha|\le m-1}
\psi_\alpha^{-1}\left(\frac{|D^{\alpha}
P_{Q_1}(x_i)-D^{\alpha}
P_{Q_2}(x_i)|}{\min\{r_1,r_2\}^{m-1-|\alpha|}}\right). \ee
\par In general the quantity $I(\cdot,\cdot)$ does not
satisfy the triangle inequality on the set
$$
\PMK:=\{T=(P,Q):~P\in\PM, Q\in\Kc\}.
$$
However, after a simple, but important modification the
function $I(\cdot,\cdot)$ transforms into a metric on
$\PMK$.
\par Namely, let us add the quantity $\rho(Q_1,Q_2)$ to the
maximum in the left-hand side of \rf{DFI}. The function
obtained we denote by $\delta$. Thus for every
$T_1=(P_1,Q_1), T_2=(P_2,Q_2)\in\PMK$ we put
\bel{D-ZM}
\delta(T_1,T_2):=\max\{\rho(Q_1,Q_2),I(P_1,P_2)\}.
\ee
\par Given $\gamma\in\R$ and $T=(P,Q)\in \PMK$ we put
$$
\gamma\circ T:=(\gamma P,Q).
$$
Now inequality \rf{SEL} is equivalent to the inequality
\bel{SDL}
\delta\left(\tfrac{1}{\lambda}\circ
T_1,\tfrac{1}{\lambda}\circ T_2\right)\le \rho(Q_1,Q_2).
\ee
(Actually, \rf{SDL} is equality, but it will be more
convenient for us to work with inequalities rather than
equalities).
\par The function $\delta$ generates  the standard {\it
geodesic metric} $d$ on $\PMK$ defined as follows: given
$T,T'\in\PMK$ we put
\bel{DFD}
d(T,T'):=\inf\sum_{i=0}^{M-1}\delta(T_i,T_{i+1})
\ee
where the infimum is taken over all finite families
$\{T_0,T_1,...,T_M\}\subset\PMK$ such that $T_0=T$ and
$T_M=T'$.
\par Our main result, Theorem \reff{T15.1}, being applied
to the case $k=0$, $\omega(t)=t^{m-1}$, states that for
every $T,T'\in\PMK$ the following inequality
$$
d(T,T')\le\delta(T,T')\le d(e^n\circ T,e^n\circ T')
$$
holds. (Of course, the first inequality is trivial and
follows from definition \rf{DFD}.) This result allows us to
reformulate Theorem \reff{ZM-LOCAP} as follows:
$f\in\ZM|_S$$\Leftrightarrow$ there exists $\lambda>0$ and
a mapping $T(Q)=(P_Q,Q)$ from $(\KS,\rho)$ into $(\PMKS,d)$
such that
\par (i) for every $Q\in\KS$ with $r_Q\le 1$ and every
$\beta, |\beta|\le m-1,$ we have
\bel{DLM} |D^{\beta}P_Q(x_Q)|r_Q^{|\beta|}\le \lambda; \ee
\par (ii) for every $Q_1,Q_2\in\KS$
\bel{DRO1} d\left(\tfrac{1}{\lambda}\circ T(Q_1),
\tfrac{1}{\lambda}\circ T(Q_2)\right)\le\rho(Q_1,Q_2); \ee
\par (iii) for every $x\in S$  we have
\bel{F-LIM}
f(x)=\lim_{x_Q=x,r_Q\to 0} P_Q(x).
\ee
\par Inequality \rf{DRO1} motivates us to introduce a
Lipschitz-type space $LO(\KS)$ of mappings
$$
T:\KS\to\Tc:=\PMKS
$$
defined by the finiteness of the following ``seminorm"
$$
\|T\|_{LO(\KS)}:=\inf\{\lambda>0:
\|\lambda^{-1}\circ T\|_{\Lip(\KS,\Tc)}\le 1\}.
$$
\par Also, inequality \rf{DLM} motivates us to define a
``norm"
$$
\|T\|^*:=\sup\{|D^{\beta}P_Q(x_Q)|r_Q^{|\beta|}:~
Q\in\KS, r_Q\le 1,|\beta|\le m-1\}.
$$
By ${\bf LO}(\KS)$ we denote a subspace of $LO(\KS)$
defined by the finiteness of the ``norm"
$$
\|T\|^*_{LO(\KS)}:=\|T\|^*+\|T\|_{LO(\KS)}.
$$
See Section 5 for details.
\par Thus one can identify the space $\ZM|_S$ with
``limiting values" $\displaystyle\lim_{x_Q=x,r_Q\to 0}
P_Q(x)$ of mappings $Q\in\KS\to T(Q)=(P_Q,Q)$ from the
space ${\bf LO}(\KS)$.
\begin{remark}\lbl{REMH}{\em By the formula:
\bel{INDKR}
\KRN\ni Q=Q(x,r)
~~\Leftrightarrow~~y=(x,r)\in\RNP \ee
we identify the family $\KRN$ of all cubes in $\RN$ with
the upper half-space $\RNP$
$$
\RNP:=\RN\times\R_+=\{y=(y_1,...,y_n,y_{n+1})\in\R^{n+1}: ~
y_{n+1}>0\}.
$$
This identification and the metric \rf{MHP} generate a
metric $\rp$ on $\RNP$ defined by the following formula:
for $z_i=(x_i,r_i)\in\RNP,i=1,2,$
\bel{DEFRP} \rp(z_1,z_2):=\ln\left(1+\frac{\max(r_1,r_2)+
\|x_1-x_2\|}{\min(r_1,r_2)}\right), \ee
if $z_1\ne z_2$, and  $\rp(z_1,z_2):=0$ whenever $z_1=z_2$.
\par Given $z=(z_1,...,z_{n+1})\in\RNP$ we put
$\bar{z}:=z=(z_1,...,-z_{n+1})$. Also, by $\|z\|_2$, we
denote the Euclidean distance in $\RNP$. We recall that the
Poincar\'e  metric on $\RNP$ is defined by the formula
\bel{DEFRH}
\rph(z_1,z_2):=\ln\frac{\|z_1-\bar{z}_2\|_2+
\|z_1-z_2\|_2}{\|z_1-\bar{z}_2\|_2- \|z_1-z_2\|_2}.
\ee
This metric is the Riemannian metric for which the line
element $ds$ is given by
$$
ds:=\frac{\sqrt{dx_1^2+...+dx_n^2+dx_{n+1}^2}}{x_{n+1}}.
$$
It determines the classical Poincar\'e upper half-space
model of the hyperbolic space $\HP:=(\RNP,\rph)$. It can be
readily seen that
\bel{RP-RH}
\rp(z_1,z_2)\sim 1+\rph(z_1,z_2),~~~
z_1,z_2\in\RNP.
\ee
\par This equivalence, \rf{DEFRP} and \rf{MHP} show that
the metric space $(\Kc,\rho)$ can be identified (up to a
constant weight) with the hyperbolic space $\HP$.
}
\end{remark}
\par In view of this remark and identification \rf{INDKR},
one can interpret the equality \rf{F-LIM} as the {\it
restriction to $\RN$} of the mapping $\RNP\ni z\to
(P_z,z)\in \PM\times \RNP$. This enables us to identify the
Zygmund space $\ZM$ with the restriction to $\RN$ of all
Lipschitz mappings of the form $T(z)=(P_z,z), z\in\RNP,$
defined on the hyperbolic space $\HP$ and taking their
values in the metric space $(\PM\times \RNP,d)$ .
\par In Section 3 and 4 we develop this approach for the
general case of the space $\CKLM$ with $k>0$. This enables
us to reformulate the extension Problem \reff{WP-CKLM} as a
geometrical problem of the existence of Lipschitz
selections of certain set-valued mappings from $\RNP$ into
a family of subsets in $\Pc_{k+m-1}\times \RNP$. We will
discuss a generalization of Problem \reff{WP-CKLM} raised
by C. Fefferman \cite{F4} (for the space $\CKO$) and show
how this problem can be reduced to the Lipschitz selection
problem for certain jet-spaces generated by functions from
$\CKLM$.
\par We observe that the Lipschitz selection method has
been used for proving the finiteness property of the spaces
$\LMT$ \cite{S1}, $\COM$ \cite{BS4} and $C^k\LMT$
\cite{BS3} (a jet-version). In \cite{S5} we used the same
technique to prove a certain weak version of the finiteness
property of the space $\CKO$.
\par {\bf Acknowledgment.} I am greatly indebted to Michael
Cwikel for helpful suggestions and remarks.
\SECT{3. The space $\CKLM$ and local polynomial
approximations}{3}
\indent \par Our notation is fairly standard. Throughout
the paper $C,C_1,C_2,...$ will be generic positive
constants which depend only on $k,m,n$. These constants can
change even in a single string of estimates. The dependence
of a constant on certain parameters is expressed, for
example, by the notation $C=C(k,m,n)$. We write $A\approx
B$ if there is a constant $C\ge 1$ such that $A/C\le B\le
CA$.
\par We let $\Pc_\ell=\Pc_\ell(\RN),$ $\ell\ge 0,$
denote the space of all polynomials on $\RN$ of degree at
most $\ell$. Finally, given $k$-times differentiable
function $f$ and $x\in\RN$, we let $T^k_x(f)$ denote the
Taylor polynomial of $f$ at $x$ of degree at most $k$:
$$
T_{x}^{k}(f)(y):=\sum_{|\alpha|\leq
k}\frac{1}{\alpha!}(D^{\alpha}f)(x)(y-x)^{\alpha}~,~~y\in
\RN.
$$
\par Finally, we put
$$
L:=k+m-1.
$$
\begin{theorem}\lbl{T14.1} Given a family of polynomials
$\{P_x\in\PK:x\in S\}$
there is a function $F\in\CKLM$ such that
$T_{x}^{k}(F)=P_{x}$ for every $x\in S$ if and only if
there is a constant $\lambda>0$ and a family of polynomials
$$
\{P_Q\in\PL:~Q\in\KS\}
$$
such that
\par (1).  For every cube
$Q\in\KS$ we have
$$
T^k_{x_Q}(P_Q)=P_{x_Q}\,;
$$
\par (2). For every $Q\in\KS$ with $r_Q\le 1$ and every
$\alpha, |\alpha|\le k,$
$$
\sup_Q|D^\alpha P_Q|\le \lambda\,;
$$
\par (3). For every $Q,Q'\in\KS$, such that $Q'\subset Q$ we
have
$$
\sup_{Q'}|P_{Q'}-P_Q|\le \lambda
(r_{Q'}+\|x_Q-x_{Q'}\|)^k\omega(r_Q).
$$
Moreover,
$$
\inf\{\|F\|_{\CKLM}:~T_{x}^{k}(F)=P_{x} ,x\in S\}\approx
\inf\lambda
$$
with constants of equivalence depending only on $k,m$ and
$n$.
\end{theorem}
\par  For the homogeneous space $\CKLM$ (normed by the
second item in \rf{14.N}) a variant of this theorem has
been proved in \cite{BS3}. The present result can be
obtained by a simple modification of the method of proof
suggested in \cite{BS3}.
\begin{theorem}\lbl{T14.2} If a function $F\in\CKLM$,
then there exists a family of polynomials
$$\{P_Q\in\PL:~Q\in\KS\}$$ such that:
\par (1).  For every cube
$Q\in\KS$ we have
$$
T^k_{x_Q}(P_Q)=T^k_{x_Q}(F);
$$
\par (2). For every $Q\in\KS$ with $r_Q\le 1$ and every
$\alpha, |\alpha|\le k,$ and $\beta, |\beta|\le
L-|\alpha|$,
$$
|D^{\alpha+\beta} P_Q(x_Q)|\le
C\|F\|_{\CKLM}\,r_Q^{-|\beta|}~;
$$
\par (3). For every
$Q_1=Q(x_1,r_1),Q_2=Q(x_2,r_2)\in\KS$, and every
$\alpha,|\alpha|\le L,$ we have
\be |D^{\alpha} P_{Q_1}(x_1)-D^{\alpha} P_{Q_2}(x_1)|&\le&
C\|F\|_{\CKLM}\nn\\&\cdot&
(\max\{r_1,r_2\}+\|x_1-x_2\|)^{L-|\alpha|}
\int\limits_{\min\{r_1,r_2\}}^{r_1+r_2+\|x_1-x_2\|}
\frac{\omega(t)}{t^m}\,dt.\nn \ee
Here  $C$ is a constant depending only on $k,m$ and $n$.
\end{theorem}
\par The proof of the theorem is based on the following
auxiliary lemmas.
\begin{lemma}\lbl{L14.3} Let $\omega\in\Omega_m$. Assume
that a family of polynomials $\{P_Q\in\PL: Q\in\KS\}$ and a
constant $\lambda>0$ satisfy the inequality
\bel{14.A}
\sup_{Q'}|P_{Q'}-P_Q|\le \lambda
(r_{Q'}+\|x_Q-x_{Q'}\|)^k\omega(r_Q),
\ee
for every $Q=Q(x_Q,r_Q),Q'=Q(x_{Q'},r_{Q'})\in\KS$ such
that $Q'\subset Q$ and $r_Q\le 4r_{Q'}$.
\par Then for every $Q,Q'\in\KS$,
$Q'\subset Q$, and every $\alpha,|\alpha|\le L,$ we have
\bel{14.C}
\sup_{Q'}|D^{\alpha} P_{Q'}-D^{\alpha} P_{Q}|\le
C\lambda \int\limits_{r_{Q'}}^{2r_Q}
\frac{\omega(t)}{t^{|\alpha|-k}}\,\frac{dt}{t}
\ee
where $C=C(k,m,n)$.
\end{lemma}
\par{\it Proof.} Let $x\in Q'$. We put
$$
Q_i:=Q(x_{Q'},2^ir_{Q'}),~i=0,...,\ell,
$$
where $\ell:=\left[\ln\frac{r_Q}{r_{Q'}}\right]+1$. Then
$$
\frac{r_{Q_\ell}}{2}\le 2^{\ell-1}r_{Q'}\le
2r_{Q}<2^{\ell}r_{Q'}=r_{Q_\ell}.
$$
Since $2r_{Q}\le r_{Q_\ell}:=2^{\ell}r_{Q'}$ and $x_{Q'}\in
Q\cap Q_\ell$, we have
\bel{14.1}
 Q\subset Q_\ell, ~~~~2r_{Q}\le r_{Q_\ell}\le 4r_{Q}.
\ee
By Markov's inequality,
$$
\sup_{Q_i}|D^{\alpha} (P_{Q_i}-P_{Q_{i+1}})|\le C\,
r_{Q_i}^{-|\alpha|}\,\sup_{Q_i}|P_{Q_i}-P_{Q_{i+1}}|\,
$$
where $C=C(k,m,n)$. Since $r_{Q_{i+1}}=2r_{Q_{i}}$ and
$x_{Q_{i+1}}=x_{Q_{i}}$, by \rf{14.A} we obtain
$$
\sup_{Q_i}|D^{\alpha} (P_{Q_i}-P_{Q_{i+1}})|\le C\lambda
r_{Q_i}^{-|\alpha|}\,(r_{Q_i}^k\omega(r_{Q_{i+1}})=
C\lambda r_{Q_i}^{k-|\alpha|}\,\omega(2r_{Q_{i}}).
$$
Since $x\in Q'\cap Q_i$, we have
$$
|D^{\alpha} P_{Q_i}(x)-D^{\alpha}P_{Q_{i+1}}(x)|\le
C\lambda
r_{Q_i}^{k-|\alpha|}\,\omega(2r_{Q_{i}}),~~~i=0,...,\ell-1.
$$
Hence
$$
|D^{\alpha} P_{Q'}(x)-D^{\alpha} P_{Q}(x)|=
|D^{\alpha} P_{Q_0}(x)-D^{\alpha}P_{Q_{\ell}}(x)|\le
\sum_{i=0}^{\ell-1}|D^{\alpha}
P_{Q_i}(x)-D^{\alpha}P_{Q_{i+1}}(x)|
$$
so that
\bel{14.2}
|D^{\alpha} P_{Q'}(x)-D^{\alpha} P_{Q}(x)|\le
C\lambda \sum_{i=0}^{\ell-1}
r_{Q_i}^{k-|\alpha|}\,\omega(2r_{Q_{i}}).
\ee
Recall that $Q\subset Q_\ell$ and $r_{Q_\ell}\le 4r_Q$, so
that by \rf{14.A} and Markov's inequality
\be |D^{\alpha} P_{Q}(x)-D^{\alpha} P_{Q_\ell}(x)| &\le&
\sup_Q |D^{\alpha} P_{Q}-D^{\alpha} P_{Q_\ell}|\nn\\
&\le& C r_Q^{-|\alpha|}\sup_Q\, \,|P_{Q}-P_{Q_\ell}|\nn\\
&\le& C\lambda (r_Q+\|x_Q-x_{Q'}\|)^k
r_Q^{-|\alpha|}\omega(r_{Q_\ell}). \nn \ee
But $x_{Q'}\in Q'\subset Q$ so that $\|x_Q-x_{Q'}\|\le
r_Q$. Therefore by \rf{14.1}
\be |D^{\alpha} P_{Q}(x)-D^{\alpha} P_{Q_\ell}(x)|&\le&
4^{|\alpha|}2^kC\lambda
r_{Q_\ell}^{k-|\alpha|}\omega(r_{Q_\ell})\nn\\
&=& 4^{|\alpha|}2^k C\lambda
(2r_{Q_{\ell-1}})^{k-|\alpha|}\omega(2r_{Q_{\ell-1}})\nn\\
&=&
C_1(k,m)\lambda
r_{Q_{\ell-1}}^{k-|\alpha|}\omega(r_{Q_{\ell-1}}).\nn \ee
Since $\omega\in\Omega_m$, we have $\omega(2t)\le
2^m\omega(t)$ so that by \rf{14.1}
\be
|D^{\alpha} P_{Q'}(x)-D^{\alpha} P_{Q}(x)|&\le&
|D^{\alpha} P_{Q'}(x)-D^{\alpha} P_{Q_\ell}(x)|+
|D^{\alpha} P_{Q_\ell}(x)-D^{\alpha} P_{Q}(x)|\nn\\
&\le& C\lambda \left(\sum_{i=0}^{\ell-1}
r_{Q_i}^{k-|\alpha|}\,\omega(2r_{Q_{i}})\right)+
C_1(k,m)\lambda
r_{Q_{\ell-1}}^{k-|\alpha|}\omega(r_{Q_{\ell-1}})
\nn\\
&\le& C_2(k,m)\lambda \sum_{i=0}^{\ell-1}
r_{Q_i}^{k-|\alpha|}\,\omega(r_{Q_{i}}).\nn \ee
Since $\omega$ is non-decreasing, for every $\alpha,
|\alpha|\ne k,$ and every $a,b, 0<a<b$, we have
$$
\int_a^b\frac{\omega(t)}{t^{|\alpha|-k}}\frac{dt}{t}\ge
\omega(a)\int_a^b\frac{1}{t^{|\alpha|-k}}\frac{dt}{t}
=(|\alpha|-k)^{-1}\omega(a) (a^{k-|\alpha|}-b^{k-|\alpha|})
$$
so that
$$
\frac{\omega(a)}{a^{|\alpha|-k}}\le \frac{(|\alpha|-k)
b^{|\alpha|-k}}{b^{|\alpha|-k}-a^{|\alpha|-k}}
\int_a^b\frac{\omega(t)}{t^{|\alpha|-k}}\frac{dt}{t}.
$$
Hence
\be r_{Q_i}^{k-|\alpha|}\,\omega(r_{Q_{i}}) &\le&
\frac{(|\alpha|-k) (2r_{Q_i})^{|\alpha|-k}}
{(2r_{Q_i})^{|\alpha|-k}-r_{Q_i}^{|\alpha|-k}}
\int_{r_{Q_i}}^{2r_{Q_i}}
\frac{\omega(t)}{t^{|\alpha|-k}}\frac{dt}{t}\nn\\ &=&
\frac{(|\alpha|-k) 2^{|\alpha|-k}} {2^{|\alpha|-k}-1}
\int_{r_{Q_i}}^{2r_{Q_i}}
\frac{\omega(t)}{t^{|\alpha|-k}}\frac{dt}{t}. \nn\ee
In a similar way we show that
$$
\omega(r_{Q_{i}}) \le \frac{1}{\ln 2}
\int_{r_{Q_i}}^{2r_{Q_i}} \omega(t)\frac{dt}{t}.
$$
Thus {\it for every} $\alpha, |\alpha|\le L$, we have
$$
r_{Q_i}^{k-|\alpha|}\,\omega(r_{Q_{i}}) \le C_3(k,m)
\int_{r_{Q_i}}^{r_{Q_{i+1}}}
\frac{\omega(t)}{t^{|\alpha|-k}}\frac{dt}{t}.
$$
\par Finally, we obtain
\be |D^{\alpha} P_{Q'}(x)-D^{\alpha} P_{Q}(x)|&\le&
C_2\lambda \sum_{i=0}^{\ell-1}
r_{Q_i}^{k-|\alpha|}\,\omega(r_{Q_{i}}) \le C_3\lambda
\sum_{i=0}^{\ell-1} \int_{r_{Q_i}}^{r_{Q_{i+1}}}
\frac{\omega(t)}{t^{|\alpha|-k}}\frac{dt}{t}\nn\\
&=& C_3\lambda \int_{r_{Q_0}}^{r_{Q_\ell}}
\frac{\omega(t)}{t^{|\alpha|-k}}\frac{dt}{t} \le C_3\lambda
\int_{r_{Q'}}^{4r_{Q}}
\frac{\omega(t)}{t^{|\alpha|-k}}\frac{dt}{t}.
\nn \ee
It can be easily seen that
$$
\int_{2r_{Q}}^{4r_{Q}}
\frac{\omega(t)}{t^{|\alpha|-k}}\frac{dt}{t} \le
C_4\int_{r_{Q}}^{2r_{Q}}
\frac{\omega(t)}{t^{|\alpha|-k}}\frac{dt}{t} \le
C_4\int_{r_{Q'}}^{2r_{Q}}
\frac{\omega(t)}{t^{|\alpha|-k}}\frac{dt}{t}
$$
so that
$$
|D^{\alpha} P_{Q'}(x)-D^{\alpha} P_{Q}(x)|\le C_3\lambda
\int_{r_{Q'}}^{4r_{Q}}
\frac{\omega(t)}{t^{|\alpha|-k}}\frac{dt}{t} \le
C_3(1+C_4)\lambda \int_{r_{Q'}}^{2r_{Q}}
\frac{\omega(t)}{t^{|\alpha|-k}}\frac{dt}{t}
$$
proving the lemma.\bx
\begin{lemma}\lbl{L14.4} Let $\lambda>0$ and let
$\{P_Q\in\PL: Q\in\KS\}$ be a family of polynomials
satisfying the conditions of Lemma \reff{L14.3}. Then for
every $Q_1=Q(x_1,r_1),Q_2=Q(x_2,r_2)\in\KS$, and every
$\alpha,|\alpha|\le L,$ we have
$$
|D^{\alpha} P_{Q_1}(x_1)-D^{\alpha} P_{Q_2}(x_1)|\le
C\lambda(\max\{r_1,r_2\}+\|x_1-x_2\|)^{L-|\alpha|}
\int\limits_{\min\{r_1,r_2\}}^{r_1+r_2+\|x_1-x_2\|}
\frac{\omega(t)}{t^m}\,dt
$$
where $C=C(k,m,n)$.
\end{lemma}
\par{\it Proof.} Put
$$
\tr:=r_1+r_2+\|x_1-x_2\|
$$
and
$$
\tQ:=Q(x_2,\tr)=Q(x_2,r_1+r_2+\|x_1-x_2\|).
$$
Then
\be
I&:=&|D^{\alpha} P_{Q_1}(x_1)-D^{\alpha}
P_{Q_2}(x_1)|\nn\\
&=& \left|\sum_{|\beta|\le p-\alpha}\frac{1}{\beta!}
[D^\beta(D^\alpha(P_{Q_2}-P_{\tQ}))](x_2)
(x_1-x_2)^\beta\right|
\nn\\
&\le& \sum_{|\beta|\le L-\alpha}\left|
D^{\alpha+\beta}(P_{Q_2}-P_{\tQ})(x_2)\right|
\|x_1-x_2\|^\beta. \nn\ee
By Lemma \reff{L14.3}
$$
|D^{\alpha+\beta} P_{Q_2}(x_2)-D^{\alpha+\beta} P_{\tQ}(x_2)| \le
C\lambda\int_{r_2}^{2\tr}
\frac{\omega(t)}{t^{|\alpha|+|\beta|-k}}\frac{dt}{t}
$$
for every $\alpha,\beta, |\alpha|+|\beta|\le L$. Hence
\be I&\le& C\lambda \sum_{|\beta|\le L-|\alpha|}
\left(\int_{r_2}^{2\tr}
\frac{\omega(t)}{t^{|\alpha|+|\beta|-k}}\frac{dt}{t}\right)
\|x_1-x_2\|^{|\beta|}\nn\\
&=& C\lambda\int_{r_2}^{2\tr} \frac{\omega(t)}{t^m}\left
(\sum_{|\beta|\le L-|\alpha|} \|x_1-x_2\|^{|\beta|}
x^{L-|\alpha|-|\beta|}\right)dt. \nn\ee
Since $\|x_1-x_2\|\le \tr=r_1+r_2+\|x_1-x_2\|$, we obtain
$$
I\le C_1\lambda \tr^{L-|\alpha|} \int_{r_2}^{2\tr}
\frac{\omega(t)}{t^m}\,dt \le C_1\lambda \tr^{L-|\alpha|}
\int_{\min\{r_1,r_2\}}^{2(r_1+r_2+\|x_1-x_2\|)}
\frac{\omega(t)}{t^m}\,dt.
$$
Since $\omega(t)/t^m$ is non-increasing, for every
$0<a<b_0\le b_1$ we have
$$
\frac{1}{b_1-a}\int_{a}^{b_1} \frac{\omega(t)}{t^m}\,dt\le
\frac{1}{b_0-a}\int_{a}^{b_0} \frac{\omega(t)}{t^m}\,dt
$$
so that
\be \int\limits_{\min\{r_1,r_2\}}^{2(r_1+r_2+\|x_1-x_2\|)}
\frac{\omega(t)}{t^m}\,dt &\le&
\frac{2(r_1+r_2+\|x_1-x_2\|)-\min\{r_1,r_2\}}
{r_1+r_2+\|x_1-x_2\|-\min\{r_1,r_2\}}
\int\limits_{\min\{r_1,r_2\}}^{r_1+r_2+\|x_1-x_2\|}
\frac{\omega(t)}{t^m}\,dt\nn\\
&=& \frac{\min\{r_1,r_2\}+2\max\{r_1,r_2\}+2\|x_1-x_2\|}
{\max\{r_1,r_2\}+\|x_1-x_2\|}
\int\limits_{\min\{r_1,r_2\}}^{r_1+r_2+\|x_1-x_2\|}
\frac{\omega(t)}{t^m}\,dt\nn\\ &\le&
3\int\limits_{\min\{r_1,r_2\}}^{r_1+r_2+\|x_1-x_2\|}
\frac{\omega(t)}{t^m}\,dt.\nn \ee
Hence
\be I&\le& 3C_1\lambda\, (r_1+r_2+\|x_1-x_2\|)^{L-|\alpha|}
\int\limits_{\min\{r_1,r_2\}}^{r_1+r_2+\|x_1-x_2\|}
\frac{\omega(t)}{t^m}\,dt\nn\\
&\le& 3\cdot 2^{L}C_1\lambda\,
(\max\{r_1,r_2\}+\|x_1-x_2\|)^{L-|\alpha|}
\int\limits_{\min\{r_1,r_2\}}^{r_1+r_2+\|x_1-x_2\|}
\frac{\omega(t)}{t^m}\,dt. \nn\ee
The lemma is proved.\bx
\medskip \par {\it Proof of Theorem \reff{T14.2}.} We put
$$
P_x:=T^k_x(F),~~x\in S.
$$
Then by Theorem \reff{T14.1} there is a family of
polynomials
$$
\{P_Q\in\PL:~Q\in\KS\}
$$
satisfying conditions (1)-(3) of this theorem with
$\lambda=C(k,m,n)\|F\|_{\CKLM}$.
\par Then equality (1) of Theorem \reff{T14.2} immediately
follows from that of Theorem \reff{T14.1}. In turn,
condition (2) of Theorem \reff{T14.1} and Markov's
inequality imply inequality (2) of Theorem \reff{T14.2}.
Finally, by Lemma \reff{L14.4}, condition (3) of Theorem
\reff{T14.1} implies inequality (3) of Theorem
\reff{T14.2}.\bx
\begin{definition}\lbl{D14.5} {\em We say that a
continuous function $\omega:\R_+\to\R_+$ is {\it
quasipower} if there is a constant $C_\omega>0$ such that
$$
\int\limits_{0}^{t} \omega(s)\frac{ds}{s}\,\le
C_\omega\,\omega(t)
$$
for all $t>0$.}
\end{definition}
\begin{example} {\em Every function $\omega(t)=t\varphi(t)$
where $\varphi$ is a non-decreasing function  is
quasi\-power (with  $C_\omega=1$). In fact
$$
\int\limits_{0}^{t} \omega(s)\frac{ds}{s}=
\int\limits_{0}^{t} \varphi(s)\,ds\le
t\varphi(t)=\omega(t).
$$
}\end{example}
\begin{theorem}\lbl{T14.7}
Let $\omega\in\Omega_m$ be a quasipower function. Suppose
that  a family of polynomials
$$\{P_Q\in\PL:~Q\in\KS\}$$
and a constant $\lambda>0$ satisfy the following
conditions:
\par (1). For every $Q\in\KS$ with $r_Q\le 1$ and every
$\alpha, |\alpha|\le k,$ and $\beta, |\beta|\le
L-|\alpha|,$ we have
$$
|D^{\alpha+\beta} P_Q(x_Q)|\le \lambda\,r_Q^{-|\beta|}~;
$$
\par (2). For every two cubes
$Q_1=Q(x_1,r_1),Q_2=Q(x_2,r_2)\in\KS$, and every
$\alpha,|\alpha|\le L,$
\bel{14.M}~~~~~~~~~ |D^{\alpha}(P_{Q_1}-P_{Q_2})(x_1)| \le
\lambda\,(\max\{r_1,r_2\}+\|x_1-x_2\|)^{L-|\alpha|}
\int\limits_{\min\{r_1,r_2\}}^{r_1+r_2+\|x_1-x_2\|}
\frac{\omega(t)}{t^m}\,dt.~~~ \ee
Then for every $x\in S$ there exists the limit
\bel{14.L}
P_x=\lim_{x_Q=x,\,r_Q\to 0}T^k_x(P_Q).
\ee
Moreover, there exists a function $F\in\CKLM$ with
$\|F\|_{\CKLM}\le C\lambda$ such that
$$
P_x=T^k_{x}(F), ~~~~x\in S,
$$
and for every $Q=Q(x,r)\in\KS$ and $|\alpha|\le k$ we have
\bel{EPF} |D^{\alpha}T^k_{x}(F)(x)-D^{\alpha}P_{Q}(x)| \le
C\lambda\,r^{k-|\alpha|}\omega(r). \ee
Here $C$ is a constant depending only on $k,m,n$ and the
constant $C_\omega$ (see Definition \reff{D14.5}).
\end{theorem}
\par {\it Proof.} By condition (2), for every $\alpha,
|\alpha|\le L,$ and every two cubes $Q_1:=Q(x,r)$ and
$Q_2:=Q(x_2,\tr)\in\KS$ with $r\le\tr\le 2r,$ we have
$$
|D^{\alpha} P_{Q_1}(x)-D^{\alpha} P_{Q_2}(x)|\le
\lambda\,\tr^{L-|\alpha|} \int\limits_{r}^{r+\tr}
\frac{\omega(t)}{t^m}\,dt\le 2^{L-|\alpha|}\lambda
r^{L-|\alpha|}\int\limits_{r}^{3r}
\frac{\omega(t)}{t^m}\,dt.
$$
Since $\omega(t)/t^m$ is non-increasing and $L:=k+m-1$, we
have
\bel{14.8} |D^{\alpha} P_{Q_1}(x)-D^{\alpha} P_{Q_2}(x)|\le
2^{L-|\alpha|}\lambda
r^{L-|\alpha|}\frac{\omega(r)}{r^m}\int\limits_{r}^{3r}
1\,dt=2^{k+m-|\alpha|}\lambda r^{k-|\alpha|}\omega(r). \ee
Consider now two cubes $Q'=Q(x,r')$ and $Q''=Q(x,r''),
r'<r''$. Put $\ell:=[\ln (r''/r')]$ and
$$
r_i:=2^ir',~~Q_i:=Q(x,r_i),~~i=0,1,...
$$
Thus
$$
r_\ell:=2^\ell r'\le r''<2^{\ell+1} r'=:r_{\ell+1}.
$$
Hence, by \rf{14.8},
\be |D^{\alpha} P_{Q'}(x)-D^{\alpha} P_{Q''}(x)|&\le&
\sum_{i=0}^{\ell-1}|D^{\alpha} P_{Q_i}(x)-D^{\alpha}
P_{Q_{i+1}}(x)|+|D^{\alpha} P_{Q_\ell}(x)-D^{\alpha}
P_{Q''}(x)|\nn\\
&\le& 2^{k+m-|\alpha|}\lambda
\left\{\left(
\sum_{i=0}^{\ell-1}r_i^{k-|\alpha|}\omega(r_i)\right)
+ r_\ell^{k-|\alpha|}\omega(r_\ell)\right\}\nn\\
&=&2^{k+m-|\alpha|}\lambda
\sum_{i=0}^{\ell}r_i^{k-|\alpha|}\omega(r_i).
\nn\ee
But
$$
\int_{r_i}^{r_{i+1}}t^{k-|\alpha|}\omega(t)\,\frac{dt}{t}\ge
\omega(r_i)\ln(r_{i+1}/r_i)=(\ln\,2)\omega(r_i),
$$
whenever $|\alpha|=k$, and
$$
\int_{r_i}^{r_{i+1}}t^{k-|\alpha|}\omega(t)\,\frac{dt}{t}\ge
\frac{\omega(r_i)}{k-|\alpha|}\left(r_{i+1}^{k-|\alpha|}
-r_{i}^{k-|\alpha|}\right)=
\frac{2^{k-|\alpha|}-1}{(k-|\alpha|)2^{k-|\alpha|}}\,
r_i^{k-|\alpha|}\omega(r_i),
$$
if  $|\alpha|\ne k$. Thus, for every $\alpha$ and every $i,
0\le i\le \ell,$ we have
$$
r_i^{k-|\alpha|}\omega(r_i)\le C(k,m)
\int_{r_i}^{r_{i+1}}t^{k-|\alpha|}\omega(t)\,\frac{dt}{t}.
$$
Hence,
\bel{14.B} |D^{\alpha} P_{Q'}(x)-D^{\alpha} P_{Q''}(x)|\le
C\lambda
\int_{r_0}^{r_{\ell}}t^{k-|\alpha|}\omega(t)\,\frac{dt}{t}
\le C\lambda
\int_{r'}^{2r''}t^{k-|\alpha|}\omega(t)\,\frac{dt}{t}. \ee
\par Consider now the case $|\alpha|\le k$. Recall that, by
the assumption, $\omega$ is a quasipower function, so that
$$
\int_{r'}^{2r''}\omega(t)\,\frac{dt}{t}\le C_\omega
\omega(2r'')\le 2^m C_\omega \omega(r'').
$$
Therefore, by \rf{14.B}, for every $\alpha, |\alpha|=k,$ we
have
$$
|D^{\alpha} P_{Q'}(x)-D^{\alpha} P_{Q''}(x)|\le C\lambda
\omega(r'')
$$
with $C=C(k,m,C_\omega)$.
\par If $|\alpha|< k$, we obtain
\be |D^{\alpha} P_{Q'}(x)-D^{\alpha} P_{Q''}(x)|&\le&
C\lambda
\int_{r'}^{2r''}t^{k-|\alpha|}\omega(t)\,\frac{dt}{t} \le
C\lambda \omega(2r'')
\int_{r'}^{2r''}t^{k-|\alpha|}\,\frac{dt}{t}\nn\\
&=& \frac{1}{k-|\alpha|}\,C\lambda \omega(2r'') \left(
(2r'')^{k-|\alpha|}-(r')^{k-|\alpha|}\right)\nn\\
&\le& \frac{1}{k-|\alpha|}\,C\lambda(2r'')^{k-|\alpha|}
\omega(2r''). \nn \ee
Thus for every $\alpha, |\alpha|\le k,$ we have
\bel{14.D} |D^{\alpha} P_{Q'}(x)-D^{\alpha} P_{Q''}(x)|\le
C\lambda(r'')^{k-|\alpha|} \omega(r''),
\ee
where $C=C(k,m,C_\omega)$.
\par Since $\omega(r)\to 0$ as $r\to 0$, there exists the
limit
\bel{14.DL}
p_\alpha(x):=\lim_{x_Q=x,\,r_Q\to 0}\,
D^{\alpha} P_{Q}(x).
\ee
We put
$$
P_x(y):=\sum_{|\beta|\le k}
\frac{p_\beta(x)}{\beta!}(y-x)^\beta.
$$
Thus $P_x\in\PK$ and
\bel{14.9}
D^\alpha P_x(x)=p_\alpha(x),~~~|\alpha|\le k.
\ee
Since $\PK$ is finite dimensional, by \rf{14.9}
$$
P_x=\lim_{x_Q=x,\,r_Q\to 0}\, T^k_xP_Q.
$$
proving \rf{14.L}.
\par Prove that the family of polynomials
$\{P_x\in\PK:~x\in S\}$ satisfies the conditions of Theorem
\reff{T14.1}, i.e., there exists a family of polynomials
$$
\{\tP_Q\in\PL:~Q\in \KS\}
$$
such that:
\par (a).  $T^k_{x_Q}(\tP_Q)=P_{x_Q}$ for every $Q\in\KS$;
\par (b). $\sup_Q|D^\alpha \tP_Q|\le C\lambda$ for every
cube $Q\in\KS$ with $r_Q\le 1$ and every $\alpha,
|\alpha|\le k$;
\par (c). For every $Q,Q'\in\KS$, such that $Q'\subset Q$ we
have
$$
\sup_{Q'}|\tP_{Q'}-\tP_Q|\le C\lambda
(r_{Q'}+\|x_Q-x_{Q'}\|)^k\omega(r_Q).
$$
\par Put
\bel{14.PW}
\tP_Q:=P_Q+P_x-T^k_xP_Q~,
~~~~Q=Q(x,r)\in\KS.
\ee
Then
\bel{14.DA} D^\alpha\tP_Q(x)=\left\{
\begin{array}{ll}
D^\alpha P_x(x),&|\alpha|\le k,\\\\
D^\alpha P_Q(x), &|\alpha|> k.
\end{array}
\right.
\ee
In particular, the condition (a) is satisfied.
\par Observe that tending $r'$ to $0$ in \rf{14.D}
we obtain
\bel{14.E}
|D^{\alpha} P_{x}(x)-D^{\alpha} P_{Q}(x)|\le
C\lambda r^{k-|\alpha|} \omega(r)
\ee
for every cube $Q=Q(x,r)\in \KS$ and every $\alpha,
|\alpha|\le k$. This proves \rf{EPF}.
\par Let us prove (b). Note that for every $\alpha,
|\alpha|\le k$, by property (1) of the theorem (with
$\beta=0$) we have
$$
|D^{\alpha} P_Q(x_Q)|\le \lambda
$$
for every $Q\in\KS$ with $r_Q\le 1$. Therefore by
\rf{14.DL} and \rf{14.9}
$$
|D^\alpha P_x(x)|=|\lim_{r\to 0}\, D^\alpha
P_{Q(x,r)}(x)|\le \lambda,~~~x\in S.
$$
Hence, by \rf{14.DA}, for every $\beta$ such that
$|\alpha+\beta|\le k$ we have
$$
|D^{\alpha+\beta} \tP_Q(x_Q)|=|D^{\alpha+\beta}
P_{x_Q}(x_Q)|\le \lambda.
$$
In turn, if $|\alpha+\beta|> k$, by \rf{14.DA} and
condition (1) of the theorem
$$
|D^{\alpha+\beta} \tP_Q(x_Q)|=|D^{\alpha+\beta}
P_{Q}(x_Q)|\le \lambda r_Q^{-|\beta|}.
$$
Therefore for every $x\in Q$ we have
\be |D^{\alpha} \tP_Q(x)|&=&|\sum_{|\beta|\le L-|\alpha|}
\frac{1}{\beta!}D^{\alpha+\beta}
\tP_Q(x_Q)(x-x_Q)^\beta|\nn\\
&\le& \sum_{|\beta|\le L-|\alpha|} |D^{\alpha+\beta}
\tP_Q(x_Q)|\,\|x-x_Q\|^{|\beta|}\nn\\
&=& \sum_{|\beta|\le k-|\alpha|} |D^{\alpha+\beta}
\tP_Q(x_Q)|\,\|x-x_Q\|^{|\beta|} \nn\\
&+& \sum_{k-|\alpha|<|\beta|\le L-|\alpha|}
|D^{\alpha+\beta}
\tP_Q(x_Q)|\,\|x-x_Q\|^{|\beta|}\nn\\
&\le& \sum_{|\beta|\le k-|\alpha|} \lambda r_Q^{|\beta|}+
\sum_{k-|\alpha|<|\beta|\le L-|\alpha|} \lambda
r_Q^{-|\beta|}r_Q^{|\beta|}. \nn \ee
Since $r_Q\le 1$, we obtain
$$
|D^{\alpha} \tP_Q(x)|\le C(k,m,n)\lambda,~~~x\in Q,
$$
proving (b).
\par Let us prove (c). Put $\br:=r+r'+\|x'-x\|$ and
$$
\bQ:=Q(x',r+r'+\|x'-x\|)=Q(x',\br).
$$
Then clearly, $Q'\subset Q\subset\bQ$, and also
$$
r\le\br=r+r'+\|x'-x\|\le r'+2r\le 3r.
$$
By \rf{14.E}
\bel{QK} |D^{\alpha} \tP_{Q'}(x')-D^{\alpha}
\tP_{\bQ}(x')|\le C\lambda \br^{k-|\alpha|} \omega(\br)\le
C\lambda r^{k-|\alpha|} \omega(r), ~~~|\alpha|\le k. \ee
On the other hand , by condition (2) of the theorem, for
every $\alpha, |\alpha|\le L$, we have
\be |D^{\alpha} P_{\bQ}(x')-D^{\alpha} P_{Q}(x')|&\le&
\lambda (\max\{\br,r\}+\|x'-x\|)^{L-|\alpha|}
\int\limits_{\min\{\br,r\}}^{\br+r+\|x'-x\|}
\frac{\omega(t)}{t^m}\,dt\nn\\
&=& \lambda (\br+\|x'-x\|)^{L-|\alpha|}
\int\limits_{r}^{\br+r+\|x'-x\|} \frac{\omega(t)}{t^m}\,dt.
\nn\ee
Since $r\le\br\le 3r,~\|x'-x\|\le r$, we obtain
\be |D^{\alpha} P_{\bQ}(x')-D^{\alpha} P_{Q}(x')|&\le&
4^{L-|\alpha|}\lambda r^{L-|\alpha|}\int\limits_{r}^{5r}
\frac{\omega(t)}{t^m}\,dt\nn\\
&\le& 4^{L-|\alpha|}\lambda
r^{L-|\alpha|}\,\,\frac{\omega(r)}{r^m}
\int\limits_{r}^{5r}1\,dt\nn\ee
so that
\bel{14.H} |D^{\alpha} P_{\bQ}(x')-D^{\alpha} P_{Q}(x')|\le
4^{k+m}\lambda r^{k-|\alpha|}\omega(r),~~~|\alpha|\le L.
\ee
Observe also that by \rf{14.E} for  all $\alpha$ with
$|\alpha|\le k$ we have
\bel{TQ}
|D^{\alpha} P_{Q}(x')-D^{\alpha}\tP_{Q}(x')|\le
C\lambda r_Q^{k-|\alpha|}\omega(r_Q)=C\lambda
r^{k-|\alpha|}\omega(r).
\ee
Hence,
\be |D^{\alpha}\tP_{Q'}(x')-D^{\alpha} \tP_{Q}(x')|&\le&
|D^{\alpha}\tP_{Q'}(x')-D^{\alpha} P_{\bQ}(x')|+
|D^{\alpha}P_{\bQ}(x')-D^{\alpha} P_{Q}(x')|\nn\\&+&
|D^{\alpha}P_{Q}(x')-D^{\alpha} \tP_{Q}(x')|\nn\ee
so that by \rf{QK}, \rf{14.H} and \rf{TQ} we obtain
\bel{14.d1}
|D^{\alpha}\tP_{Q'}(x')-D^{\alpha}
\tP_{Q}(x')|\le C\lambda r^{k-|\alpha|}\omega(r)~~~{\rm for
~all}~~\alpha,|\alpha|\le k.
\ee
\par Now consider the case $k<|\alpha|\le L$. By
\rf{14.DA} and \rf{14.B} we have
\be
|D^{\alpha}\tP_{Q'}(x')-D^{\alpha} \tP_{\bQ}(x')|&=&
|D^{\alpha}P_{Q'}(x')-D^{\alpha} P_{\bQ}(x')|\nn\\ &\le&
C\lambda \int\limits_{r'}^{2\br}
\frac{\omega(t)}{t^{|\alpha|-k}}\,\frac{dt}{t} \le C\lambda
\omega(2\br)\int\limits_{r'}^{2\br}
t^{k-|\alpha|}\,\frac{dt}{t}\nn\\&\le&
C\lambda\,(|\alpha|-k)^{-1} (r')^{k-|\alpha|}\omega(2\br).
\nn\ee
Since $\br\le 3r$, we obtain
$$
|D^{\alpha}\tP_{Q'}(x')-D^{\alpha} \tP_{\bQ}(x')|\le
C\lambda\,\frac{2^m}{|\alpha|-k}
(r')^{k-|\alpha|}\omega(r)=C_1\lambda\,
(r')^{k-|\alpha|}\omega(r).
$$
On the other hand, by \rf{14.H} and \rf{14.DA}
$$
|D^{\alpha}\tP_{\bQ}(x')-D^{\alpha}\tP_{Q}(x')|=
|D^{\alpha}P_{\bQ}(x')-D^{\alpha}P_{Q}(x')|\le
 C\lambda\,
r^{k-|\alpha|}\omega(r).
$$
Hence
\be |D^{\alpha}\tP_{Q'}(x')-D^{\alpha} \tP_{Q}(x')|&\le&
|D^{\alpha}\tP_{Q'}(x')-D^{\alpha} \tP_{\bQ}(x')|+
|D^{\alpha}\tP_{\bQ}(x')-D^{\alpha}\tP_{Q}(x')|\nn\\
&\le& C\lambda\, ((r')^{k-|\alpha|}\omega(r)
+r^{k-|\alpha|}\omega(r))\nn\\
&\le& C\lambda\,
(r')^{k-|\alpha|}\omega(r)),~~~k<|\alpha|\le L. \nn\ee
\par We have proved that for every two cubes $Q'=Q(x',r')$
and $Q=Q(x,r)$ such that $Q'\subset Q$, we have
\bel{14.FD} |D^{\alpha}\tP_{Q'}(x')-D^{\alpha} \tP_{Q}(x')|
\le C\lambda\,\left\{
\begin{array}{ll}
r^{k-|\alpha|}\omega(r),&|\alpha|\le k,\\\\
(r')^{k-|\alpha|}\omega(r)), &k<|\alpha|\le L.
\end{array}
\right.
\ee
\par To estimate $\sup_{Q'}|\tP_{Q'}-\tP_{Q}|$
we let $\tQ$ denote the cube $\tQ:=Q(x,r'+\|x-x'\|)$. Since
$Q'\subset Q$, we have $r'+\|x-x'\|\le r$ so that
$Q'\subset \tQ\subset Q$. Also, by \rf{14.DA},
$$
D^{\alpha}\tP_{\tQ}(x)=D^{\alpha}P_{x}(x)
=D^{\alpha}\tP_{Q}(x)
$$
so that by \rf{14.FD}
$$ |D^{\alpha}\tP_{\tQ}(x')-D^{\alpha}
\tP_{Q}(x)| \le C\lambda\,\left\{
\begin{array}{ll}
0,&|\alpha|\le k,\\\\
(r'+\|x-x'\|)^{k-|\alpha|}\omega(r), &k<|\alpha|\le L.
\end{array}
\right.
$$
Therefore for every $y\in \tQ$ we have
\be |\tP_{\tQ}(y)-\tP_{Q}(y)|&=& \left|\sum_{|\alpha|\le
L}\frac{1}{\alpha!} \left(D^\alpha\tP_{\tQ}(x)-
D^\alpha\tP_{Q}(x)\right)(y-x)^\alpha\right|\nn\\
&=& \left|\sum_{k<|\alpha|\le L}\frac{1}{\alpha!}
\left(D^\alpha\tP_{\tQ}(x)-
D^\alpha\tP_{Q}(x)\right)(y-x)^\alpha\right|\nn\\
&\le& \sum_{k<|\alpha|\le L} \left|D^\alpha\tP_{\tQ}(x)-
D^\alpha\tP_{Q}(x)\right|\|y-x\|^{|\alpha|}\nn\\
&\le& C\lambda\, \sum_{k<|\alpha|\le L}
(r'+\|x-x'\|)^{k-|\alpha|}\omega(r)\,r_{\tQ}^{|\alpha|}\nn\\
&\le&
C_1\lambda\, (r'+\|x-x'\|)^k\omega(r). \nn\ee
Thus
$$
\sup_{Q'}|\tP_{\tQ}-\tP_{Q}| \le
\sup_{\tQ}|\tP_{\tQ}-\tP_{Q}|\le C\lambda\,
(r'+\|x-x'\|)^k\omega(r).
$$
Let us estimate $\sup_{Q'}|\tP_{Q'}-\tP_{\tQ}|$. By
\rf{14.FD} for every $\alpha, |\alpha|\le L,$
$$ |D^{\alpha}\tP_{Q'}(x')-D^{\alpha}
\tP_{\tQ}(x)| \le C\lambda\,\left\{
\begin{array}{ll}
r_{\tQ}^{k-|\alpha|}\omega(r_{\tQ}),&|\alpha|\le k,\\\\
(r')^{k-|\alpha|}\omega(r_{\tQ}), &k<|\alpha|\le L.
\end{array}
\right.
$$
Hence, for each $y\in Q'$ we have
\be |\tP_{Q'}(y)-\tP_{\tQ}(y)|&=& \left|\sum_{|\alpha|\le
L}\frac{1}{\alpha!} \left(D^\alpha\tP_{Q'}(x')-
D^\alpha\tP_{\tQ}(x')\right)(y-x')^\alpha\right|\nn\\
&\le& \sum_{|\alpha|\le L} \left|D^\alpha\tP_{Q'}(x')-
D^\alpha\tP_{\tQ}(x')\right|\|y-x'\|^{|\alpha|}\nn\\
&=& \sum_{|\alpha|\le k} \left|D^\alpha\tP_{Q'}(x')-
D^\alpha\tP_{\tQ}(x')\right|\|y-x'\|^{|\alpha|}\nn\\&+&
\sum_{k<|\alpha|\le L} \left|D^\alpha\tP_{Q'}(x')-
D^\alpha\tP_{\tQ}(x')\right|\|y-x'\|^{|\alpha|}\nn\\
&\le& C\lambda\, \sum_{|\alpha|\le k}
(r_{\tQ})^{k-|\alpha|}\omega(r_{\tQ})(r')^{k-|\alpha|}\nn\\
&+& C\lambda\, \sum_{k<|\alpha|\le L}
(r')^{k-|\alpha|}\omega(r_{\tQ})(r')^{|\alpha|}. \nn\ee
Since $r'\le r_{\tQ}=r'+\|x-x'\|$, we obtain
$$
|\tP_{Q'}(y)-\tP_{\tQ}(y)|\le C\lambda\,
r_{\tQ}^{k}\omega(r_{\tQ})= C\lambda\,
(r'+\|x-x'\|)^{k}\omega(r_{\tQ})
$$
proving that
$$
\sup_{Q'}|\tP_{Q'}-\tP_{\tQ}|\le C\lambda\,
(r'+\|x-x'\|)^{k}\omega(r_{\tQ}).
$$
Finally, since $r_{\tQ}=r'+\|x-x'\|\le r$, we have
\be \sup_{Q'}|\tP_{Q'}-\tP_{Q}|&\le&
\sup_{Q'}|\tP_{Q'}-\tP_{\tQ}|+
\sup_{Q'}|\tP_{\tQ}-\tP_{Q}|\nn\\ &\le& C\lambda\,
(r'+\|x-x'\|)^{k}\omega(r_{\tQ})+ C\lambda\,
(r'+\|x-x'\|)^{k}\omega(r)\nn\\
&\le& 2C\lambda\,
(r'+\|x-x'\|)^{k}\omega(r). \nn\ee
Theorem \ref{T14.7} is completely proved.\bx
\SECT{4. $\CKLM$ as a space of Lipschitz mappings}{4}
\indent \par The point of departure for our approach is the
inequality \rf{14.M} of Theorem \reff{T14.7}. This
inequality motivates the definition of a certain metric on
the set
$$
\PL\times\KRN=\{T=(P,Q):~ P\in\PL,~Q\in\KRN\}.
$$
\par This allows us to identify the restriction $\CKLM|_S$
with a space of Lipschitz mappings from $\KS$ (equipped
with a certain hyperbolic-type metric) into $\PKR$.
\par Given $v>0$ and a multiindex $\alpha, |\alpha|\le L,$
we define a function $\fa(\cdot;v)$ on $\R_+$ by letting
\bel{15.0} \fa(t;v):=t^{L-|\alpha|}\int_v^{v+t}
\frac{\omega(s)}{s^m}\,ds. \ee
(Recall that $L:=k+m-1$.) By $\fin(\cdot;v)$ we denote the
inverse to the function $\fa(\cdot;v)$ (i.e., the inverse
to the function $\fa$ with respect to the first argument).
Since for every $v>0$ the function $\fa(\cdot;v)$ is
strictly increasing, the function $\fin(\cdot;v)$ is
well-defined.
\par Thus for every $u\ge 0$ we have
\bel{15.1} \fin(u;v)^{L-|\alpha|}
\int\limits_v^{v+\fin(u;v)}\frac{\omega(s)}{s^m}\,ds=u. \ee
In particular,
\bel{15.AP}
\int\limits_v^{v+\fin(u;v)}\frac{\omega(s)}{s^m}\,ds=u,
~~~|\alpha|=L. \ee
\par Now fix two elements
$$
T_1=(P_1,Q_1),~T_2=(P_2,Q_2)\in\PKR
$$
where $Q_1=Q(x_1,r_1),Q_2=Q(x_2,r_2)\in\KRN$ and
$P_i\in\PL,i=1,2$. Put
\be \Delta(T_1,T_2)&:=& \max
\{\max\{r_1,r_2\}+\|x_1-x_2\|,\label{DDELTA}\\
&&\max_{|\alpha|\le L,\,i=1,2}
\fin(|D^\alpha(P_1-P_2)(x_i)|;\min\{r_1,r_2\})\}
 \nn\ee
and
\bel{15.2} \domg(T_1,T_2):= \int\limits_{\min\{r_1,r_2\}}^
{\min\{r_1,r_2\}+\Delta(T_1,T_2)}
\frac{\omega(s)}{s^m}\,ds~, \ee
if $T_1\ne T_2$, and $\domg(T_1,T_2):=0$, if $T_1=T_2$.
\par Observe that definition \rf{DDELTA} and equality
\rf{15.AP} imply the following explicit formula for
$\domg(T_1,T_2),~ T_1\ne T_2,$:
\be \domg(T_1,T_2)&:=& \max\left\{
\int\limits_{\min\{r_1,r_2\}}^{r_1+r_2+\|x_1-x_2\|}
\frac{\omega(s)}{s^m}\,ds,~
\max_{|\alpha|=L}|D^\alpha P_1-D^\alpha P_2|,\right.\nn\\
&& \left.\max_{|\alpha|<
L,\,i=1,2}\int\limits_{\min\{r_1,r_2\}}^
{\min\{r_1,r_2\}+\fin(|D^\alpha(P_1-P_2)(x_i)|;
\min\{r_1,r_2\})} \frac{\omega(s)}{s^m}\,ds\right\}. \nn\ee
\par Let us introduce a function
$\rom:\KRN\times\KRN\to\R_+$ by letting
\bel{DROQ} \rom(Q_1,Q_2):=\left\{
\begin{array}{ll}
\int\limits_{\min\{r_1,r_2\}}^{r_1+r_2+\|x_1-x_2\|}
\frac{\omega(s)}{s^m}\,ds\,, & Q_1\ne Q_2,\\\\
0,& Q_1=Q_2.
\end{array}
\right.
\ee
Here $Q_i=Q(x_i,r_i),~i=1,2.$
As we shall see below, $\rom$ is a metric on $\KRN$, see
Remark \reff{REMRO}.
\par Observe that for every $P\in\PL$ and
$Q_i\in\KRN,i=1,2,$ we have
\bel{15.EP}
\domg((P,Q_1),(P,Q_2))=\rom(Q_1,Q_2).
\ee
\par In these settings the inequality \rf{14.M} of Theorem
\reff{T14.7} can be reformulated in the following way.
\begin{claim}\lbl{CLAIM} Given a family of polynomials
$$\{P_Q\in\PL:~Q\in\KS\}$$
and a constant $\lambda>0$ the following two statements are
equivalent:
\par (i). For every two cubes
$Q_1=Q(x_1,r_1),Q_2=Q(x_2,r_2)\in\KS$, and every
$\alpha,|\alpha|\le L,$
\bel{IN14.7}~~~~~~~ |(D^{\alpha} P_{Q_1}-D^{\alpha}
P_{Q_2})(x_1)|\le \lambda\,(\max\{r_1,r_2\}+\|x_1-x_2\|)^L
\int\limits_{\min\{r_1,r_2\}}^{r_1+r_2+\|x_1-x_2\|}
\frac{\omega(t)}{t^m}\,dt. \ee
(Recall that this is inequality \rf{14.M} of Theorem
\reff{T14.7}).
\par (ii). Let $T:\KS\to\PKR$ be a mapping defined by the
formula $T(Q):=(P_Q,Q),~Q\in\KS$. Then
\bel{DRO}
\domg(\lambda^{-1}\circ T(Q_1),\lambda^{-1}\circ
T(Q_2))\le \rom(Q_1,Q_2),~~Q_1,Q_2\in\KS.
\ee
\end{claim}
\par (Recall that $\lambda\circ T:=(\lambda P,Q)$ provided
$T=(P,Q)\in\PKR$ and $\lambda\in\R$.)
\par {\it Proof.} Put
$$
A_i:=|(D^{\alpha} (\lambda^{-1}P_{Q_1})-D^{\alpha}
(\lambda^{-1}P_{Q_2}))(x_i)|,~~~i=1,2.
$$
By definition \rf{15.1} inequality \rf{IN14.7} can be
reformulated as follows:
$$
A_1\le \fa(\max\{r_1,r_2\}+\|x_1-x_2\|;\min\{r_1,r_2\}).
$$
Hence
$$
\fin(A_1;\min\{r_1,r_2\})\le \max\{r_1,r_2\}+\|x_1-x_2\|.
$$
Changing the order of cubes in this inequality and taking
the maximum over all $\alpha,|\alpha|\le L,$ we conclude
that \rf{IN14.7} is equivalent to the following inequality:
$$
\max_{|\alpha|\le L,\,i=1,2}\fin(A_i;\min\{r_1,r_2\})\le
\max\{r_1,r_2\}+\|x_1-x_2\|.
$$
In turn, by definition \rf{DDELTA}, this inequality is
equivalent to the next one:
\be \Delta\left(\tfrac{1}{\lambda}\circ
T(Q_1),\tfrac{1}{\lambda}\circ
T(Q_2)\right)&=&\max\limits_{|\alpha|\le L,\,i=1,2}
\{\max\{r_1,r_2\}+\|x_1-x_2\|,
\fin(A_i;\min\{r_1,r_2\})\}\nn\\
&\le&\max\{r_1,r_2\}+\|x_1-x_2\|.
\nn \ee
\par Since the function
$t\to\int_v^{v+t}\frac{\omega(s)}{s^m}\,ds$ is strictly
increasing, the inequality
$$
\Delta(\lambda^{-1}\circ T(Q_1),\lambda^{-1}\circ
T(Q_2))\le\max\{r_1,r_2\}+\|x_1-x_2\|
$$
is equivalent to
\be
\domg(\lambda^{-1}\circ T(Q_1),\lambda^{-1}\circ
T(Q_2)) &:=& \int\limits_{\min\{r_1,r_2\}}^
{\min\{r_1,r_2\}+\Delta(\lambda^{-1}\circ
T(Q_1),\lambda^{-1}\circ T(Q_2))} \frac{\omega(s)}{s^m}\,ds
\nn\\
&\le& \int\limits_{\min\{r_1,r_2\}}^
{\min\{r_1,r_2\}+\max\{r_1,r_2\}+\|x_1-x_2\|}
\frac{\omega(s)}{s^m}\,ds\nn\\
&=&\int\limits_{\min\{r_1,r_2\}}^ {r_1+r_2+\|x_1-x_2\|}
\frac{\omega(s)}{s^m}\,ds = \rom(Q_1,Q_2).
\nn\ee
The claim is proved.\bx
\par Let us define a metric on $\PKR$ as a geodesic metric
generated by the function $\domg$. Given $T,T'\in\PKR$ we
put
\bel{dm} \dm(T,T'):=\inf\sum_{i=0}^{M-1}\dom(T_i,T_{i+1})
\ee
where the infimum is taken over all finite families
$\{T_0,T_1,...,T_M\}\subset\PKR$ such that $T_0=T$ and
$T_M=T'$.
\par  Observe several elementary properties of $\dm$. In
particular, as we have noted above, the function
$\rom:\KRN\times\KRN\to\R_+$ is a {\it metric} on $\KRN$,
see Remark \reff{REMRO}. This property of $\rom$ and
\rf{15.2} and \rf{DDELTA} immediately imply the following
inequality: for every $P_i\in\PL, Q_i\in\KRN, i=1,2,$ we
have
$$
\dm((P_1,Q_1),(P_2,Q_2))\ge \rom(Q_1,Q_2).
$$
In turn, this inequality and \rf{15.EP} imply the
following:
$$
\dm((P,Q_1),(P,Q_2))=\rom(Q_1,Q_2),~~~P\in\PL,
Q_i\in\KRN,~i=1,2.
$$
\par The main result of the section is the following
\begin{theorem}\lbl{T15.1} For every $T,T'\in\PKR$ we have
$$
\dm(T,T')\le \dom(T,T') \le \dm(e^n\circ T,e^n\circ T').
$$
\end{theorem}
The proof of this theorem relies on a series of auxiliary
lemmas.
\begin{lemma}\lbl{L15.2} For every
$b_0,b_1,...,b_\ell>0$, $a_0,a_1,...,a_{\ell-1}\ge 0$ and
$c_0,c_1,...,c_{\ell-1}\ge 0$ we have
\be
&\max&\left\{\int\limits_{\min\{b_0,b_\ell\}}^{b_0+b_\ell+
\sum\limits_{i=0}^{\ell-1}a_i}\frac{\omega(t)}{t^m}\,dt,
\int\limits_{\min\{b_0,b_\ell\}}^{\min\{b_0,b_\ell\}+
\sum\limits_{i=0}^{\ell-1}c_i}
\frac{\omega(t)}{t^m}\,dt\right\}\nn\\
&\le&\sum\limits_{i=0}^{\ell-1} \max\left\{
\int\limits_{\min\{b_i,b_{i+1}\}}^{b_i+b_{i+1}+a_i}
\frac{\omega(t)}{t^m}\,dt,
\int\limits_{\min\{b_i,b_{i+1}\}}^{\min\{b_i,b_{i+1}\}+c_i}
\frac{\omega(t)}{t^m}\,dt\right\}. \nn\ee
\end{lemma}
\par {\it Proof.} We put $s_{-1}=a_{-1}=c_{-1}:=0,$
$$
s_i:=\max\{\max\{b_i,b_{i+1}\}+a_i,c_i\},~~~i=0,...,\ell-1,
$$
and
$$
I:=[\min\{b_0,b_\ell\},\min\{b_0,b_\ell\}+
\max\{\max\{b_0,b_\ell\}+\sum\limits_{i=0}^{\ell-1}a_i,
\sum\limits_{i=0}^{\ell-1}c_i\}].
$$
Then the inequality of the lemma is equivalent to the
following one:
\bel{LFI}
\int\limits_I\frac{\omega(t)}{t^m}\,dt \le
\sum\limits_{i=0}^{\ell-1}
\int\limits_{\min\{b_i,b_{i+1}\}}^{\min\{b_i,b_{i+1}\}+s_i}
\frac{\omega(t)}{t^m}\,dt.
\ee
\par To prove this inequality we put
$$
I_i:=[\min\{b_0,b_\ell\}+\sum\limits_{j=-1}^{i-1}s_j,
\min\{b_0,b_\ell\}+\sum\limits_{j=-1}^{i}s_j],~~~i=0,...,
\ell-1.
$$
Then $|I_i|=s_i$ and
\be
\sum\limits_{j=-1}^{l-1}s_j&=& \sum\limits_{j=-1}^{l-1}
\max\{\max\{b_j,b_{j+1}\}+a_j,c_j\}\nn\\
&\ge&
\max\{\sum\limits_{j=-1}^{l-1}\max\{b_j,b_{j+1}\}+
\sum\limits_{j=-1}^{l-1}a_j,
\sum\limits_{j=-1}^{l-1}c_j\}\nn\\
&\ge& \max\{\max\{b_0,b_{\ell}\}+
\sum\limits_{j=-1}^{l-1}a_j,\sum\limits_{j=-1}^{l-1}c_j\}.
\nn\ee
Thus
$$
\bigcup\limits_{i=0}^{\ell-1}I_i\supset I
$$
so that
\bel{LI1}
\int\limits_I\frac{\omega(t)}{t^m}\,dt \le
\sum\limits_{i=0}^{\ell-1} \int\limits_{I_i}
\frac{\omega(t)}{t^m}\,dt.
\ee
\par Put
$$
A_i:=[\min\{b_i,b_{i+1}\},\min\{b_i,b_{i+1}\}+s_i],~~~
i=0,...,\ell-1.
$$
Then $|A_i|=|I_i|=s_i$. But the left end of the segment
$I_i$ is bigger than the left end of the segment $A_i$. In
fact,
$$
\min\{b_0,b_{\ell}\}+\sum\limits_{j=-1}^{i-1}s_j\ge s_{i-1}
=\max\{\max\{b_{i-1},b_i\}+a_{i-1},c_{i-1}\}\ge b_i\ge
\min\{b_i,b_{i+1}\}.
$$
Thus the segment $A_i$ is a shift of $I_i$ to the left.
Since $\omega(t)/t^m$ is non-increasing, this implies
$$
\int\limits_{I_i}\frac{\omega(t)}{t^m}\,dt \le
\int\limits_{A_i} \frac{\omega(t)}{t^m}\,dt=
\int\limits_{\min\{b_i,b_{i+1}\}}^{\min\{b_i,b_{i+1}\}+s_i}
\frac{\omega(t)}{t^m}\,dt.
$$
This inequality and inequality \rf{LI1} imply \rf{LFI}. The
lemma is proved.\bx
\begin{remark}\lbl{REMRO} {\em We put in Lemma \reff{L15.2}
$\ell=2, c_0=c_1=c_2:=0$ and get
\bel{LRO} \int\limits_{\min\{b_0,b_2\}}^{b_0+b_2+
a_0+a_1}\frac{\omega(t)}{t^m}\,dt \le
\int\limits_{\min\{b_0,b_1\}}^{b_0+b_1+a_0}
\frac{\omega(t)}{t^m}\,dt +
\int\limits_{\min\{b_1,b_2\}}^{b_1+b_2+a_1}
\frac{\omega(t)}{t^m}\,dt. \ee
This inequality easily implies the triangle inequality for
the function $\rom:\KRN\times\KRN\to\R_+$ defined by
\rf{DROQ}. In fact, for every cubes $Q_i=Q(x_i,r_i)\in\KRN,
i=0,1,2,$ we have
$$
\rom(Q_0,Q_2)\le \int\limits_{\min\{r_0,r_1\}}
^{r_0+r_2+\|x_0-x_2\|} \frac{\omega(s)}{s^m}\,ds \le
\int\limits_{\min\{r_0,r_2\}}^
{r_0+r_2+\|x_0-x_1\|+\|x_1-x_2\|} \frac{\omega(s)}{s^m}\,ds
$$
so that by \rf{LRO}
$$
\rom(Q_0,Q_2)\le\int\limits_{\min\{r_0,r_1\}}
^{r_0+r_1+\|x_0-x_1\|} \frac{\omega(s)}{s^m}\,ds+
\int\limits_{\min\{r_1,r_2\}}^ {r_1+r_2+\|x_1-x_2\|}
\frac{\omega(s)}{s^m}\,ds= \rom(Q_0,Q_1)+\rom(Q_1,Q_2).
$$
}\end{remark}
\begin{lemma}\lbl{L15.4} Let $\{P_0,P_1,...,P_\ell\}$ be a
finite subfamily of $\PL$ and let $\{x_0,x_1,...,x_\ell\}$
be a subset of $\RN$.
\par Then for every $\alpha, |\alpha|\le L,$ we have
$$
|D^\alpha(P_0-P_\ell)(x_0)|\le e^n\max_{|\beta|\le
L-|\alpha|} \left(\sum\limits_{i=0}^{\ell-1}
|D^{\alpha+\beta}(P_i-P_{i+1})(x_i)|\right)\cdot
\left(\sum\limits_{i=0}^{\ell-1}
\|x_i-x_{i+1}\|\right)^{|\beta|}.
$$
\end{lemma}
\par {\it Proof.} We have
\be |D^{\alpha}(P_i-P_{i+1})(x')| &=& |\sum\limits_{|\beta|
\le L-|\alpha|}\frac{1}{\beta!}D^{\alpha+\beta}
(P_i-P_{i+1})(x_i)\cdot (x'-x_i)^\beta|\nn\\
&\le& \sum\limits_{|\beta| \le
L-|\alpha|}\frac{1}{\beta!}|D^{\alpha+\beta}
(P_i-P_{i+1})(x_i)| \cdot\|x'-x_i\|^{|\beta|} \nn \ee
so that
\be |D^{\alpha}(P_0-P_{\ell})(x')| &\le&
\sum\limits_{i=0}^{\ell-1}|D^{\alpha}(P_i-P_{i+1})(x')|\nn\\
&\le& \sum\limits_{i=0}^{\ell-1}\,\, \sum\limits_{|\beta|
\le L-|\alpha|}\frac{1}{\beta!}|D^{\alpha+\beta}
(P_i-P_{i+1})(x_i)| \cdot\|x'-x_i\|^{|\beta|}\nn\\
&=& \sum\limits_{|\beta| \le
L-|\alpha|}\,\,\frac{1}{\beta!} \left(
\sum\limits_{i=0}^{\ell-1}|D^{\alpha+\beta}
(P_i-P_{i+1})(x_i)| \cdot\|x'-x_i\|^{|\beta|} \right). \nn
\ee
Hence,
\be |D^{\alpha}(P_0-P_{\ell})(x')| &\le& \left(
\sum\limits_{|\beta| \le L}\frac{1}{\beta!}\right) \max
\limits_{|\beta| \le L-|\alpha|}
\sum\limits_{i=0}^{\ell-1}|D^{\alpha+\beta}
(P_i-P_{i+1})(x_i)| \cdot\|x'-x_i\|^{|\beta|}\nn\\
&\le& e^n \max \limits_{|\beta| \le L-|\alpha|}
\sum\limits_{i=0}^{\ell-1}|D^{\alpha+\beta}
(P_i-P_{i+1})(x_i)| \cdot\|x'-x_i\|^{|\beta|}. \nn \ee
It remains to note that
$$
\|x'-x_i\|=\|x_0-x_i\|\le \sum\limits_{j=0}^{\ell-1}
\|x_j-x_{j+1}\|
$$
and the lemma follows.\bx
\begin{lemma}\lbl{L15.5} For every $v,R,t>0$ and every
$\alpha,\beta, |\alpha+\beta|\le L,$ we have
$$
\int\limits_{v}^ {v+\fin(R^{|\beta|}t;v)}
\frac{\omega(s)}{s^m}\,ds \le \max\left\{\int\limits_{v}^
{v+R} \frac{\omega(s)}{s^m}\,ds\, , \int\limits_{v}^
{v+\varphi^{-1}_{\alpha+\beta}(t;v)}
\frac{\omega(s)}{s^m}\,ds\right\}\,.
$$
\end{lemma}
\par {\it Proof.} Put $s:=\varphi^{-1}_{\alpha+\beta}(t;v)$
and $q:=\max\{R,s\}$. Since $R\le q$, we obtain
$$
R^{|\beta|}q^{L-|\alpha+\beta|}\le q^{L-|\alpha|},
$$
so that
$$
R^{|\beta|}q^{L-|\alpha+\beta|}
\int\limits_{v}^{v+q}\frac{\omega(s)}{s^m}\,ds \le
q^{L-|\alpha|}
\int\limits_{v}^{v+q}\frac{\omega(s)}{s^m}\,ds \,.
$$
By definition, see \rf{15.0},
$$
\fa(u;v):=
u^{L-|\alpha|}\int_v^{v+u}\frac{\omega(s)}{s^m}\,ds
\,,~~~u>0,
$$
so that the latter inequality can be written in the
following form:
$$
R^{|\beta|}\varphi_{\alpha+\beta}(q;v)\le \fa(q;v).
$$
Since $s\le q$ and $\varphi_{\alpha+\beta}$ is increasing,
this inequality implies the following
$$
R^{|\beta|}\varphi_{\alpha+\beta}(s;v)\le \fa(q;v).
$$
But $\varphi_{\alpha+\beta}(s;v)=t$ so that
$R^{|\beta|}t\le \fa(q;v)$. Since $\fin$ is increasing, we
have
$$
\fin(R^{|\beta|}t;v)\le
q=\max\{R,\varphi^{-1}_{\alpha+\beta}(t;v)\}
$$
proving the lemma.\bx
\begin{lemma}\lbl{L15.MFI} For every
$b_0,b_1,...,b_\ell>0$, $u_0,u_1,...,u_{\ell-1}\ge 0$ and
every $\alpha, |\alpha|\le L,$ we have
\be
&&\int\limits_{\min\{b_0,b_\ell\}}^{\min\{b_0,b_\ell\}+
\fin(\sum\limits_{i=0}^{\ell-1}u_i; \min\{b_0,b_\ell\})}
\frac{\omega(t)}{t^m}\,dt \le\nn\\
&&\sum\limits_{i=0}^{\ell-1}
\max\left\{\int\limits_{\min\{b_i,b_{i+1}\}}^{b_i+b_{i+1}}
\frac{\omega(t)}{t^m}\,dt,
\int\limits_{\min\{b_i,b_{i+1}\}}^{\min\{b_i,b_{i+1}\}+
\fin(u_i; \min\{b_i,b_{i+1}\})}
\frac{\omega(t)}{t^m}\,dt\right\}.
\nn
\ee
\end{lemma}
\par {\it Proof.} Put
$$
A:=\fin\left(\sum\limits_{i=0}^{\ell-1}u_i;
\min\{b_0,b_\ell\}\right)
$$
and
$$
c_i:=\fin(u_i;\min\{b_i,b_{i+1}\}),~~~i=0,...,\ell-1.
$$
\par Assume that
$$
A\le
B:=\sum\limits_{i=0}^{\ell-1}c_i=
\sum\limits_{i=0}^{\ell-1}\fin(u_i;\min\{b_i,b_{i+1}\}).
$$
Hence
\be
I&:=&\int\limits_{\min\{b_0,b_\ell\}}^{\min\{b_0,b_\ell\}+
\fin(\sum\limits_{i=0}^{\ell-1}u_i; \min\{b_0,b_\ell\})}
\frac{\omega(t)}{t^m}\,dt\,
=\int\limits_{\min\{b_0,b_\ell\}}^{\min\{b_0,b_\ell\}+A}
\frac{\omega(t)}{t^m}\,dt\nn\\
&\le&
\int\limits_{\min\{b_0,b_\ell\}}^{\min\{b_0,b_\ell\}+B}
\frac{\omega(t)}{t^m}\,dt\,=
\int\limits_{\min\{b_0,b_\ell\}}^{\min\{b_0,b_\ell\}+
\sum\limits_{i=0}^{\ell-1}c_i} \frac{\omega(t)}{t^m}\,dt
\nn\ee
so that by Lemma \reff{L15.2} (with $a_i=0,
i=0,...,\ell-1$) we obtain
$$
I\le \sum\limits_{i=0}^{\ell-1}
\max\left\{\int\limits_{\min\{b_i,b_{i+1}\}}^{b_i+b_{i+1}}
\frac{\omega(t)}{t^m}\,dt\,,
\int\limits_{\min\{b_i,b_{i+1}\}}^{\min\{b_i,b_{i+1}\}+c_i}
\frac{\omega(t)}{t^m}\,dt\right\}.
$$
This proves the lemma under the assumption $A\le B$.
\par Suppose that $A>B$. By identity \rf{15.1}
$$
\left(\fin\left(\sum\limits_{i=0}^{\ell-1}u_i;
\min\{b_0,b_\ell\}\right)\right)^{L-|\alpha|}
\int\limits_{\min\{b_0,b_\ell\}}^{\min\{b_0,b_\ell\}+
\fin\left(\sum\limits_{i=0}^{\ell-1}u_i;
\min\{b_0,b_\ell\}\right)}\frac{\omega(t)}{t^m}\,dt =
\sum\limits_{i=0}^{\ell-1}u_i
$$
so that
$$
A^{L-|\alpha|}
\int\limits_{\min\{b_0,b_\ell\}}^{\min\{b_0,b_\ell\}+A}
\frac{\omega(t)}{t^m}\,dt=\sum\limits_{i=0}^{\ell-1}u_i.
$$
Hence
$$
I=\int\limits_{\min\{b_0,b_\ell\}}^{\min\{b_0,b_\ell\}+A}
\frac{\omega(t)}{t^m}\,dt=A^{|\alpha|-L}
\left(A^{L-|\alpha|}
\int\limits_{\min\{b_0,b_\ell\}}^{\min\{b_0,b_\ell\}+A}
\frac{\omega(t)}{t^m}\,dt\right)=
A^{|\alpha|-L}\sum\limits_{i=0}^{\ell-1}u_i.
$$
Since $A>B$ and $|\alpha|\le L$, we obtain
\bel{EI} I=A^{|\alpha|-L}\sum\limits_{i=0}^{\ell-1}u_i \le
B^{|\alpha|-L}\sum\limits_{i=0}^{\ell-1}u_i . \ee
\par Again, by identity \rf{15.1} for every
$i=0,...,\ell-1,$ we have
$$
u_i=(\fin(u_i;\min\{b_i,b_{i+1}\}))^{L-|\alpha|}
\int\limits_{\min\{b_i,b_{i+1}\}}^{\min\{b_i,b_{i+1}\}+
\fin(u_i;\min\{b_i,b_{i+1}\})}\frac{\omega(t)}{t^m}\,dt
$$
so that by \rf{EI}
$$
I\le \sum\limits_{i=0}^{\ell-1}
\left(\frac{\fin(u_i;\min\{b_i,b_{i+1}\})}{B}\right)
^{L-|\alpha|}
\int\limits_{\min\{b_i,b_{i+1}\}}^{\min\{b_i,b_{i+1}\}+
\fin(u_i;\min\{b_i,b_{i+1}\})}\frac{\omega(t)}{t^m}\,dt\,.
$$
But
$$
\fin(u_i;\min\{b_i,b_{i+1}\})\le B:=
\sum\limits_{j=0}^{\ell-1}\fin(u_j;\min\{b_j,b_{j+1}\})
$$
for every $i=0,...,\ell-1,$ so that
$$
I\le \sum\limits_{i=0}^{\ell-1}
\int\limits_{\min\{b_i,b_{i+1}\}}^{\min\{b_i,b_{i+1}\}+
\fin(u_i;\min\{b_i,b_{i+1}\})}\frac{\omega(t)}{t^m}\,dt\,.
$$
The lemma is proved.\bx
\medskip
\par {\it Proof of Theorem \reff{T14.1}}. The inequality
$\dm(T,T')\le\domg(T,T')$ trivially follows from definition
\rf{dm} of the metric $\dm$.
\par In turn, the inequality $\domg(T,T')\le\dm(e^n\circ
T,e^n\circ T')$ is equivalent to the following statement:
Let
$$
\{T_i=(P_i,Q_i)\in\PKR:~i=0,...,\ell\}
$$
where $Q_i=Q(x_i,r_i)$, be a subfamily of $\PKR$  such that
$T_0=T, T_\ell=T'$. Put
\bel{DI}
I:=\int\limits_{\min\{r_0,r_\ell\}}^{\min\{r_0,r_\ell\}+
\Delta(T_0,T_\ell)}\frac{\omega(t)}{t^m}\,dt\,,
\ee
and
$$
A:=\sum\limits_{i=0}^{\ell-1}
\int\limits_{\min\{r_i,r_{i+1}\}}^{\min\{r_i,r_{i+1}\}+
\Delta(e^n\circ T_i,e^n\circ T_{i+1})}
\frac{\omega(t)}{t^m}\,dt\,.
$$
Then
\bel{IA}
I\le A.
\ee
For the sake of brevity we put
$$
v:=\min\{r_0,r_\ell\}.
$$
Recall that
\be
\Delta(T_0,T_\ell)&:=&\max_{|\alpha|\le L}
\{\max\{r_0,r_\ell\}+\|x_0-x_\ell\|,\nn\\
&&\fin(|D^\alpha(P_0-P_\ell)(x_0)|;v),
\fin(|D^\alpha(P_0-P_\ell)(x_\ell)|;v)\}.
\nn\ee
Given multiindex $\alpha, |\alpha|\le L,$ we put
$$
I_\alpha:= \int\limits_{v}^{v+
\fin(|D^\alpha(P_0-P_\ell)(x_0)|;v)}
\frac{\omega(t)}{t^m}\,dt\,,
$$
and
$$
J_\alpha:= \int\limits_{v}^{v+
\fin(|D^\alpha(P_0-P_\ell)(x_\ell)|;v)}
\frac{\omega(t)}{t^m}\,dt\,.
$$
Then by \rf{DI}
\bel{FINQ} I:=\int\limits_{v}^{v+\Delta(T_0,T_\ell)}
\frac{\omega(t)}{t^m}\,dt =\max_{|\alpha|\le L}
\left\{\int\limits_{v}^{v+
\|x_0-x_\ell\|}\frac{\omega(t)}{t^m}\,dt\,, I_\alpha,
J_\alpha\right\}. \ee
\par Prove that
\bel{E1} \int\limits_{v}^{v+
\|x_0-x_\ell\|}\frac{\omega(t)}{t^m}\,dt\,\le A. \ee
In fact,
$$
\int\limits_{v}^{v+
\sum\limits_{i=0}^{\ell-1}\|x_i-x_{i+1}\|}
\frac{\omega(t)}{t^m}\,dt=
\int\limits_{\min\{r_0,r_\ell\}}^{\min\{r_0,r_\ell\}+
\sum\limits_{i=0}^{\ell-1}\|x_i-x_{i+1}\|}
\frac{\omega(t)}{t^m}\,dt\le
\int\limits_{\min\{r_0,r_\ell\}}^{r_0+r_\ell+
\sum\limits_{i=0}^{\ell-1}\|x_i-x_{i+1}\|}
\frac{\omega(t)}{t^m}\,dt
$$
so that applying Lemma \reff{L15.2} (with $c_i=0,
i=0,...,\ell-1,$) we obtain
$$
\int\limits_{v}^{v+
\sum\limits_{i=0}^{\ell-1}\|x_i-x_{i+1}\|}
\frac{\omega(t)}{t^m}\,dt\le \sum\limits_{i=0}^{\ell-1}
\int\limits_{\min\{r_i,r_{i+1}\}}
^{r_i+r_{i+1}+\|x_i-x_{i+1}\|} \frac{\omega(t)}{t^m}\,dt.
$$
By definition \rf{DDELTA} for every $i=0,...,\ell-1,$ we
have
$$
\max\{r_i,r_{i+1}\}+\|x_i-x_{i+1}\|\le \Delta(e^n\circ
T_i,e^n\circ T_{i+1})
$$
so that
\be r_i+r_{i+1}+\|x_i-x_{i+1}\|&=& \min\{r_i,r_{i+1}\}+
\max\{r_i,r_{i+1}\}+\|x_i-x_{i+1}\|\nn\\ &\le&
\min\{r_i,r_{i+1}\}+\Delta(e^n\circ T_i,e^n\circ T_{i+1}).
\nn\ee
Hence
\bel{INTI} \int\limits_{\min\{r_i,r_{i+1}\}}
^{r_i+r_{i+1}+\|x_i-x_{i+1}\|} \frac{\omega(t)}{t^m}\,dt
\le \int\limits_{\min\{r_i,r_{i+1}\}}
^{\min\{r_i,r_{i+1}\}+\Delta(e^n\circ T_i,e^n\circ
T_{i+1})} \frac{\omega(t)}{t^m}\,dt \ee
which implies the following inequality
\bel{EI1} \int\limits_{v}^{v+
\sum\limits_{i=0}^{\ell-1}\|x_i-x_{i+1}\|}
\frac{\omega(t)}{t^m}\,dt\le \sum\limits_{i=0}^{\ell-1}
\int\limits_{\min\{r_i,r_{i+1}\}}^{\min\{r_i,r_{i+1}\}+
\Delta(e^n\circ T_i,e^n\circ T_{i+1})}
\frac{\omega(t)}{t^m}\,dt\,=A. \ee
Thus
$$
\int\limits_{v}^{v+
\|x_0-x_\ell\|}\frac{\omega(t)}{t^m}\,dt\,\le
\int\limits_{v}^{v+
\sum\limits_{i=0}^{\ell-1}\|x_i-x_{i+1}\|}
\frac{\omega(t)}{t^m}\,dt\le A.
$$
proving \rf{E1}.
\par Now prove that
\bel{E2} I_\alpha\le A,~~~|\alpha|\le L. \ee
To this end given multiindex $\gamma$ and $i=0,...,\ell-1,$
we put
\bel{DG}
U_{\gamma,i}:=|D^\gamma(e^nP_i-e^nP_{i+1})(x_i)|
\ee
and
$$
U_{\gamma}:=\sum\limits_{i=0}^{\ell-1}U_{\gamma,i}.
$$
We also set
$$
R:=\sum\limits_{i=0}^{\ell-1}\|x_i-x_{i+1}\|.
$$
Then by Lemma \reff{L15.4}
\be |D^\alpha(P_0-P_{\ell})(x_0)| &\le& \max_{|\beta|\le
L-|\alpha|} \left(\sum\limits_{i=0}^{\ell-1}
|D^{\alpha+\beta}(e^nP_i-e^nP_{i+1})(x_i)|\right)\cdot
\left(\sum\limits_{i=0}^{\ell-1}
\|x_i-x_{i+1}\|\right)^{|\beta|}\nn\\
&=& \max_{|\beta|\le L-|\alpha|} R^{|\beta|}
U_{\alpha+\beta}\,. \nn\ee
Hence
$$
I_\alpha\le \int\limits_{v}^{v+
\fin\left(\max\limits_{|\beta|\le L-|\alpha|}R^{|\beta|}
U_{\alpha+\beta};v\right)} \frac{\omega(t)}{t^m}\,dt=
\max_{|\beta|\le L-|\alpha|}\int\limits_{v}^{v+
\fin(R^{|\beta|}
U_{\alpha+\beta};v)}\frac{\omega(t)}{t^m}\,dt
$$
so that by Lemma \reff{L15.5}
\bel{IAL} I_\alpha\le \max_{|\beta|\le L-|\alpha|}\left\{
\int\limits_{v}^{v+R}\frac{\omega(t)}{t^m}\,dt\,,
\int\limits_{v}^{v+ \varphi^{-1}_{\alpha+\beta}
(U_{\alpha+\beta};v)}\frac{\omega(t)}{t^m}\,dt\right\}. \ee
By \rf{EI1}
\bel{LA} \int\limits_{v}^{v+R}\frac{\omega(t)}{t^m}\,dt\,=
\int\limits_{v}^{v+\sum\limits_{i=0}^{\ell-1}
\|x_i-x_{i+1}\|} \frac{\omega(t)}{t^m}\,dt\,\le A. \ee
\par Prove that
\bel{SE} \int\limits_{v}^{v+ \varphi^{-1}_{\alpha+\beta}
(U_{\alpha+\beta};v)}\frac{\omega(t)}{t^m}\,dt\le A. \ee
For the sake of brevity we put
$$
v_i:=\min\{r_i,r_{i+1}\},~~~i=0,...,\ell-1.
$$
Then by Lemma \reff{L15.MFI}
\be \int\limits_{v}^{v+ \varphi^{-1}_{\alpha+\beta}
(U_{\alpha+\beta};v)}\frac{\omega(t)}{t^m}\,dt &=&
\int\limits_{\min\{r_0,r_\ell\}}^{\min\{r_0,r_\ell\}+
\varphi^{-1}_{\alpha+\beta}
(\sum\limits_{i=0}^{\ell-1}U_{\alpha+\beta,i};
\min\{r_0,r_\ell\})} \frac{\omega(t)}{t^m}\,dt\nn\\
&\le& \sum\limits_{i=0}^{\ell-1}\max\left\{
\int\limits_{v_i} ^{r_i+r_{i+1}}
\frac{\omega(t)}{t^m}\,dt\,, \int\limits_{v_i}
^{v_i+\varphi^{-1}_{\alpha+\beta}(U_{\alpha+\beta,i};v_i)}
\frac{\omega(t)}{t^m}\,dt \right\}. \nn \ee
By \rf{INTI}
\bel{R1} \int\limits_{v_i}^{r_i+r_{i+1}}
\frac{\omega(t)}{t^m}\,dt=
\int\limits_{\min\{r_i,r_{i+1}\}}^{r_i+r_{i+1}}
\frac{\omega(t)}{t^m}\,dt \le
\int\limits_{\min\{r_i,r_{i+1}\}}
^{\min\{r_i,r_{i+1}\}+\Delta(e^n\circ T_i,e^n\circ
T_{i+1})} \frac{\omega(t)}{t^m}\,dt. \ee
In turn, by \rf{DG} and definition \rf{DDELTA}
$$
\varphi^{-1}_{\alpha+\beta}(U_{\alpha+\beta,i};v_i)=
\varphi^{-1}_{\alpha+\beta}
(|D^\gamma(e^nP_i-e^nP_{i+1})(x_i)|;\min\{r_i,r_{i+1}\})
\le \Delta(e^n\circ T_i,e^n\circ T_{i+1})
$$
so that
$$
\int\limits_{v_i}^{v_i+\varphi^{-1}_{\alpha+\beta}
(U_{\alpha+\beta,i};v_i)}
\frac{\omega(t)}{t^m}\,dt \le
\int\limits_{\min\{r_i,r_{i+1}\}}
^{\min\{r_i,r_{i+1}\}+\Delta(e^n\circ T_i,e^n\circ
T_{i+1})} \frac{\omega(t)}{t^m}\,dt.
$$
Hence
$$
\int\limits_{v}^{v+ \varphi^{-1}_{\alpha+\beta}
(U_{\alpha+\beta};v)}\frac{\omega(t)}{t^m}\,dt \le
\sum\limits_{i=0}^{\ell-1}\int\limits_{\min\{r_i,r_{i+1}\}}
^{\min\{r_i,r_{i+1}\}+\Delta(e^n\circ T_i,e^n\circ
T_{i+1})} \frac{\omega(t)}{t^m}\,dt=A,
$$
proving \rf{SE}.
\par Now inequality \rf{E2} follows from
\rf{IAL}, \rf{LA} and \rf{SE}. In the same way we prove
that
$$
J_\alpha\le A,~~~|\alpha|\le L.
$$
Finally, this inequality, \rf{SE},\rf{E1} and \rf{FINQ}
imply the required inequality \rf{IA}.
\par Theorem \reff{T14.1} is proved.\bx
\medskip
\par We present several results related to calculation of
the function $\domg$, see \rf{DDELTA} and \rf{15.2}, and
the metric $\dm$, see \rf{dm}. \par Let
$T_1=(P_1,Q_1),~T_2=(P_2,Q_2)\in\PKR,$ where
$Q_1=Q(x_1,r_1),Q_2=Q(x_2,r_2)\in\KRN$ and
$P_i\in\PL,i=1,2$. Fix a point $y\in\RN$ and put
\be
\Delta(T_1,T_2;y)&:=& \max
\{\max\{r_1,r_2\}+\|x_1-x_2\|,\label{DEFDY}\\
&&\max_{|\alpha|\le L}
\fin(|D^\alpha(P_1-P_2)(y)|;\min\{r_1,r_2\})\}. \nn\ee
Thus, we define $\Delta(T_1,T_2;y)$ by replacing in
\rf{DDELTA} the points $x_i, i=1,2,$ by $y$ .
\par We also put
\bel{DELY} \domg(T_1,T_2;y):=
\int\limits_{\min\{r_1,r_2\}}^
{\min\{r_1,r_2\}+\Delta(T_1,T_2;y)}
\frac{\omega(s)}{s^m}\,ds~, \ee
for $T_1\ne T_2$, and $\domg(T_1,T_2;y):=0$, whenever
$T_1=T_2$. Clearly,
$$
\Delta(T_1,T_2)=\max\{\Delta(T_1,T_2;x_1),
\Delta(T_1,T_2;x_2)\},
$$
and
$$
\domg(T_1,T_2)=\max\{\domg(T_1,T_2;x_1),
\domg(T_1,T_2;x_2)\}.
$$
\begin{proposition}\lbl{PDN} For every $y,z\in\RN$ and
every $T_i=(P_i,Q_i)\in\PKR,$ where
$Q_i=Q(x_i,r_i)\in\KRN,~i=1,2,$ we have
$$
\domg(T_1,T_2;z)\le\domg(\gamma\circ T_1,\gamma\circ
T_2;y),
$$
where
$$
\gamma=\max\left\{1,
\frac{e^n\|y-z\|^L}{(\max\{r_1,r_2\}+\|x_1-x_2\|)^L}\right
\}.
$$
\end{proposition}
\par {\it Proof.} Fix a multiindex $\alpha, |\alpha|\le L,$
and put
$$
\tP_0:=P_1,~\tP_1=P_1,~\tP_2:=P_2,~~~~{\rm and}~~~~
\tx_0:=z,~\tx_1=y,~\tx_2:=y.
$$
Let us apply Lemma \reff{L15.4} to polynomials
$\{\tP_0,\tP_1,\tP_2\}$ and points $\{\tx_0,\tx_1,\tx_2\}.$
We have
\be
|D^\alpha(P_1-P_2)(z)|&=&|D^\alpha(\tP_0-\tP_2)(\tx_0)|
\nn\\
&\le& e^n\max_{|\beta|\le L-|\alpha|}
\left(\sum\limits_{i=0}^{1}
|D^{\alpha+\beta}(\tP_i-\tP_{i+1})(\tx_i)|\right)\cdot
\left(\sum\limits_{i=0}^{1}
\|\tx_i-\tx_{i+1}\|\right)^{|\beta|}\nn\\
&=& e^n\max_{|\beta|\le L-|\alpha|}
|D^{\alpha+\beta}(P_1-P_2)(y)|\cdot \|z-y\|^{|\beta|}.
\nn
\ee
Put
$$
U_{\alpha+\beta}:=|D^{\alpha+\beta}(\gamma P_1-\gamma
P_2)(y)|.
$$
Then
\be |D^\alpha(P_1-P_2)(z)| &\le& \max_{|\beta|\le
L-|\alpha|} |D^{\alpha+\beta}(P_1-P_2)(y)|\nn\\
&\cdot&
e^n\left(\frac{\|z-y\|}{\max\{r_1,r_2\}+\|x_1-x_2\|}
\right)^{|\beta|}(\max\{r_1,r_2\}+\|x_1-x_2\|)
^{|\beta|}\nn\\
&\le& \max_{|\beta|\le L-|\alpha|} |D^{\alpha+\beta}(\gamma
P_1-\gamma P_2)(y)|\cdot (\max\{r_1,r_2\}+\|x_1-x_2\|)
^{|\beta|}\nn\\
&=& \max_{|\beta|\le L-|\alpha|}
(\max\{r_1,r_2\}+\|x_1-x_2\|) ^{|\beta|} U_{\alpha+\beta}.
\nn
\ee
We also put
$$
v:=\min\{r_1,r_2\},~~~A:=\max\{r_1,r_2\}.
$$
Then
$$
I_\alpha:= \int\limits_{v}^
{v+\fin(|D^\alpha(P_1-P_2)(z)|;v)}
\frac{\omega(s)}{s^m}\,ds \le \max_{|\beta|\le L-|\alpha|}
\int\limits_{v}^ {v+\fin((A+\|x_1-x_2\|) ^{|\beta|}
U_{\alpha+\beta};v)} \frac{\omega(s)}{s^m}\,ds,
$$
so that by Lemma \reff{L15.5}
$$
I_\alpha\le \max_{|\beta|\le L-|\alpha|}\left\{
\int\limits_{v}^ {v+A+\|x_1-x_2\|}
\frac{\omega(s)}{s^m}\,ds, \int\limits_{v}^
{v+\varphi^{-1}_{\alpha+\beta}(U_{\alpha+\beta};v)}
\frac{\omega(s)}{s^m}\,ds\right\}.
$$
But by \rf{DEFDY}
$$
A+\|x_1-x_2\|=\max\{r_1,r_2\}+\|x_1-x_2\|\le
\Delta(\gamma\circ T_1,\gamma\circ T_2;y),
$$
and
$$
\varphi^{-1}_{\alpha+\beta}(U_{\alpha+\beta};v)=
\varphi^{-1}_{\alpha+\beta}(|D^{\alpha+\beta}(\gamma
P_1-\gamma P_2)(y)|;\min\{r_1,r_2\}) \le \Delta(\gamma\circ
T_1,\gamma\circ T_2;y).
$$
Hence
$$
I_\alpha\le
 \int\limits_{v}^
{v+\Delta(\gamma\circ T_1,\gamma\circ T_2;y)}
\frac{\omega(s)}{s^m}\,ds=\domg(\gamma\circ T_1,\gamma\circ
T_2;y),
$$
proving the proposition.\bx
\par Proposition \reff{PDN} and Theorem \reff{T15.1} imply
the following
\begin{corollary}\lbl{CYZ} Let $\theta\ge 1$. For
every $T_i=(P_i,Q_i)\in\PKR,$ where
$Q_i=Q(x_i,r_i)\in\KRN,~i=1,2,$ and every point $y\in\RN$
such that
$$
\|y-x_1\|\le \theta\,(r_1+r_2+\|x_1-x_2\|)
$$
we have
$$
\domg(\gamma^{-1}\circ T_1,\gamma^{-1}\circ T_2;y)
\le
\domg(T_1,T_2)\le\domg(\gamma\circ T_1,\gamma\circ T_2;y),
$$
and
$$
\domg(\gamma^{-1}\circ T_1,\gamma^{-1}\circ T_2;y) \le
\dm(T_1,T_2)\le\domg(\gamma\circ T_1,\gamma\circ T_2;y).
$$
Here $\gamma=\gamma(n,\theta)$ is a constant depending only
on $n$ and $\theta$.
\end{corollary}
\par For instance, by this corollary,
\bel{DX12}
\domg(\gamma^{-1}\circ T_1,\gamma^{-1}\circ
T_2;x_i) \le \dm(T_1,T_2)\le\domg(\gamma\circ
T_1,\gamma\circ T_2;x_i),~~~i=1,2,
\ee
or
$$
\domg\left(\gamma^{-1}\circ T_1,\gamma^{-1}\circ
T_2;\frac{x_1+x_2}{2}\right) \le
\dm(T_1,T_2)\le\domg\left(\gamma\circ T_1,\gamma\circ
T_2;\frac{x_1+x_2}{2}\right),
$$
with $\gamma=\gamma(n)$ depending only on $n$.
\par Let us present one more formula for calculation of
$\domg$. Given $v>0$ and multiindex $\alpha, |\alpha|\le
L,$ we put
$$
h(t;v):=\int\limits_{v}^{v+t}\frac{\omega(s)}{s^m}\,ds\,,
~~~t>0,
$$
and
$$
\psi_\alpha(u;v):=\left(t\left[h^{-1}(t;v)
\right]^{L-|\alpha|}\right)^{-1}(u),~~~u>0.
$$
Thus
\bel{PME} t\left[h^{-1}(t;v)\right]^{L-|\alpha|}
=\psi_\alpha^{-1}(t;v),~~~t>0. \ee
\begin{proposition} For every $T_i=(P_i,Q_i)\in\PKR,$ where
$Q_i=Q(x_i,r_i)\in\KRN,~i=1,2,$ we have
$$
\domg(T_1,T_2;x_1)=\max_{|\alpha|\le L}
\left\{\int\limits_{\min\{r_1,r_2\}}^
{r_1+r_2+\|x_1-x_2\|}\frac{\omega(s)}{s^m}\,ds\,,
\psi_\alpha(|D^\alpha(P_1-P_2)(x_1)|;\min\{r_1,r_2\})
\right\}.
$$
\end{proposition}
\par {\it Proof.} We put $v:=\min\{r_1,r_2\}$,
$$
U_\alpha:=|D^\alpha(P_1-P_2)(x_1)|,
$$
and
$$
I_\alpha:=\int\limits_{v}^
{v+\fin(U_\alpha;v)}\frac{\omega(s)}{s^m}\,ds\,.
$$
Recall that
$$
\domg(T_1,T_2;x_1):=\max_{|\alpha|\le L}
\left\{\int\limits_{v}^
{r_1+r_2+\|x_1-x_2\|}\frac{\omega(s)}{s^m}\,ds\,,
I_\alpha\right\}.
$$
Prove that
\bel{PSU} I_\alpha=\psi_\alpha(U_\alpha;v). \ee
Recall that by \rf{15.1}
\bel{DEM}
\fin(u;v)^{L-|\alpha|}
\int\limits_v^{v+\fin(u;v)}\frac{\omega(s)}{s^m}\,ds=u
\ee
so that
$$
I_\alpha=\frac{\fin(U_\alpha;v)^{L-|\alpha|}}
{\fin(U_\alpha;v)^{L-|\alpha|}}\int\limits_{v}^
{v+\fin(U_\alpha;v)}\frac{\omega(s)}{s^m}\,ds\,
=\frac{U_\alpha}{\fin(U_\alpha;v)^{L-|\alpha|}}.
$$
Hence
$$
\fin(U_\alpha;v)= \left(\frac{U_\alpha}{I_\alpha}\right)^
{\frac{1}{L-|\alpha|}}
$$
so that by \rf{DEM}
$$
\left(\left(\frac{U_\alpha}{I_\alpha}\right)^
{\frac{1}{L-|\alpha|}}\right)^{L-|\alpha|}
\int\limits_v^{v+\left(\frac{U_\alpha}{I_\alpha}\right)^
{\frac{1}{L-|\alpha|}}}\frac{\omega(s)}{s^m}\,ds
=U_\alpha.
$$
We obtain
$$
\frac{U_\alpha}{I_\alpha}
\int\limits_v^{v+\left(\frac{U_\alpha}{I_\alpha}\right)^
{\frac{1}{L-|\alpha|}}}\frac{\omega(s)}{s^m}\,ds =U_\alpha
$$
which implies
$$
I_\alpha=
\int\limits_v^{v+\left(\frac{U_\alpha}{I_\alpha}\right)^
{\frac{1}{L-|\alpha|}}}\frac{\omega(s)}{s^m}\,ds =
h((U_\alpha/I_\alpha)^ {\frac{1}{L-|\alpha|}};v).
$$
Hence $ h^{-1}(I_\alpha;v)=(U_\alpha/I_\alpha)^
{\frac{1}{L-|\alpha|}},$ so that
$I_\alpha\,[h^{-1}(I_\alpha;v)]^{L-|\alpha|}= U_\alpha.$ By
\rf{PME} this implies
$\psi_\alpha^{-1}(I_\alpha;v)=U_\alpha$ proving \rf{PSU}
and the proposition.\bx
\par In particular, by this proposition and \rf{DX12},
for every $T_i=(P_i,Q_i)\in\PKR,$ where
$Q_i=Q(x_i,r_i)\in\KRN,~i=1,2,$ we have
\be \dm\left(\tfrac{1}{\gamma}\circ
T_1,\tfrac{1}{\gamma}\circ T_2\right)&\le&
\max_{|\alpha|\le L}\left\{\int\limits_{\min\{r_1,r_2\}}^
{r_1+r_2+\|x_1-x_2\|}\frac{\omega(s)}{s^m}\,ds,
\psi_\alpha(|D^\alpha(P_1-P_2)(x_1)|;\min\{r_1,r_2\})
\right\}\nn\\
&\le& \dm(\gamma\circ T_1,\gamma\circ T_2) \nn\ee
where $\gamma=\gamma(n)$ depends only on $n$.
\par Let us present two examples.
\begin{example} {\em Consider the case of the Zygmund space
$\ZM:=\Lambda^m_\omega(\RN)$ with $\omega(t)=t^{m-1}$, see
\rf{NR-ZM}. In this case
$$
h(t;v):=\int\limits_{v}^{v+t}\frac{1}{s}\,ds=
\ln\left(1+\frac{t}{v}\right), ~~~t>0,
$$
so that $h^{-1}(s;v)=v(e^s-1)$. Hence
$$
\psi_\alpha(u;v):=
\left[sv^{m-1-|\alpha|}(e^s-1)^{m-1-|\alpha|}
\right]^{-1}(u)= \left[s(e^s-1)^{m-1-|\alpha|}\right]^{-1}
\left(\frac{u}{v^{m-1-|\alpha|}}\right).
$$
Thus in this case for every $T_i=(P_i,Q_i)\in\PKR,$ where
$Q_i=Q(x_i,r_i)\in\KRN,~i=1,2,$ we have
\be
\domg(T_1,T_2)&=&\max_{|\alpha|\le m-1,i=1,2}
\left\{\ln\left(1+\frac{\max\{r_1,r_2\}+\|x_1-x_2\|}
{\min\{r_1,r_2\}}\right),\right.\nn\\
&& \left.\left[s(e^s-1)^{m-1-|\alpha|}
\right]^{-1}\left(\frac{|D^\alpha(P_1-P_2)(x_i)|}
{\min\{r_1,r_2\}^{m-1-|\alpha|}}\right)\right\}, \nn\ee
}\end{example}
so that we again obtain the formula \rf{D-ZM} for
$\delta_\omega$.
\begin{example} \lbl{EX-COM}{\em Consider the space
$\CKO=C^k\Lambda^1_\omega(\RN)$ with $\omega(t)=t$. This is
the Sobolev space $W^{k+1}_\infty(\RN)$ of bounded
functions $f\in\CK$ whose partial derivatives of order $k$
satisfy the Lipschitz condition:
$$
|D^{\alpha}f(x)-D^{\alpha}f(y)|\le\lambda\|x-y\|,
~~x,y\in\RN,~~|\alpha|=k.
$$
Since $m=1$ and $L=k$, we have
$$
h(t;v):=\int\limits_{v}^{v+t}1\,ds=t, ~~~t>0,
$$
and
$$
\psi_\alpha(u;v):=\left[ss^{k-|\alpha|}\right]^{-1}(u)=
u^{\frac{1}{k+1-|\alpha|}}.
$$
Thus, in this case for every $T_i=(P_i,Q_i)\in\PK\times
\Kc,$ where $Q_i=Q(x_i,r_i)\in\KRN,~i=1,2,$ we have
$$
\domg(T_1,T_2)=\max_{|\alpha|\le m-1,i=1,2}
\left\{\max\{r_1,r_2\}+\|x_1-x_2\|,
|D^\alpha(P_1-P_2)(x_i)|^{\frac{1}{m-|\alpha|}} \right\}.
$$
\par A function $\domg$ of such a kind and the metric $\dm$
generated by $\domg$ have been studied in \cite{S5}.
}\end{example}
\par Now by Claim \reff{CLAIM} and Theorem \reff{T15.1},
the results of Theorem \reff{T14.2} and Theorem
\reff{T14.7} can be formulated in the following form.
\begin{theorem}\lbl{T14.2N} Let $F\in\CKLM$ and let
$\lambda:=\|F\|_{\CKLM}$. There exists a family of
polynomials $\{P_Q\in\PL:~Q\in\KS\}$ such that:
\par (1).  $T^k_{x_Q}(P_Q)=T^k_{x_Q}(F)$ for every cube
$Q\in\KS$;
\par (2). For every $Q\in\KS$ with $r_Q\le 1$ and every
$\alpha, |\alpha|\le k,$ and $\beta, |\beta|\le L$,
$$
|D^{\alpha+\beta} P_Q(x_Q)|\le C\lambda\,r_Q^{-|\beta|}~;
$$
\par (3). The mapping $T(Q):=(P_Q,Q),~Q\in\KS,$
satisfies the Lipschitz condition
\bel{LipT}
\dm((C\lambda)^{-1}\circ
T(Q_1),(C\lambda)^{-1}\circ T(Q_2))\le
\rom(Q_1,Q_2),~~Q_1,Q_2\in\KS.
\ee
Here  $C$ is a constant depending only on $k,m$ and $n$.
\end{theorem}
\begin{theorem}\lbl{T14.7N}
Let $\omega\in\Omega_m$ be a quasipower function. Assume
that  a mapping
$$
T(Q)=(P_Q,Q),~~~~Q\in\KS,
$$
from $\KS$ into $\PKR$ and a constant $\lambda>0$ satisfy
the following conditions:
\par (1). For every $Q\in\KS$ with $r_Q\le 1$ and every
$\alpha, |\alpha|\le k,$ and $\beta, |\beta|\le
L-|\alpha|,$ we have
$$
|D^{\alpha+\beta} P_Q(x_Q)|\le \lambda\,r_Q^{-|\beta|}~;
$$
\par (2). For every $Q_1,Q_2\in\KS$
\bel{BL}
\dm(\lambda^{-1}\circ T(Q_1),\lambda^{-1}\circ
T(Q_2))\le \rom(Q_1,Q_2).
\ee
Then there exists a function $F\in\CKLM$ with
$\|F\|_{\CKLM}\le C\lambda$ such that for every $x\in S$
$$
T^k_{x}(F)=\lim_{x_Q=x,\,r_Q\to 0}T^k_x(P_Q).
$$
Moreover, for every $Q=Q(x,r)\in\KS$ and every $\alpha$,
$|\alpha|\le k$ we have
$$
|D^{\alpha}T^k_{x}(F)(x)-D^{\alpha}P_{Q}(x)| \le
C\lambda\,r^{k-|\alpha|}\omega(r).
$$
\par Here the constant $C$ depends only on $k,m,n$ and the
constant $C_\omega$. \end{theorem}
\par Recall that the metric $\rom$ is defined by formula
\rf{DROQ} and the constant $C_\omega$ is defined in
Definition \reff{D14.5}.
\par Theorem \reff{T14.2N} and Theorem \reff{T14.7N} enable
us to describe the space $\CKLM$ and its restrictions to
subsets of $\RN$ as a certain space of Lipschitz mappings
defined on subsets of the metric space
$$
\KRN_\omega:=(\KRN,\rom)
$$
and taking their values in the metric space
$$
\TOM:=(\PKR,\dm).
$$
\par Let $\MRR$ be a metric space and let $\Lip(\Mcc,\TOM)$
be the space of Lipschitz mappings from $\Mcc$ into $\PKR$
equipped with the standard Lipschitz seminorm
$$
\|T\|_{\Lip(\Mcc,\TOM)}:=\sup_{x,y\in\Mcc,\,x\ne y}
\frac{\dm(T(x),T(y))}{\rho(x,y)}.
$$
We introduce a Lipschitz-type space $LO(\Mcc,\TOM)$ of
mappings $T:\Mcc\to\PKR$ defined by the finiteness of the
following seminorm:
\bel{DLO}
\|T\|_{LO(\Mcc,\TOM)}:=\inf\{\lambda>0:~\|\lambda^{-1}\circ
T\|_{\Lip(\Mcc,\TOM)}\le 1\}. \ee
\par Thus $T\in LO(\Mcc,\TOM)$ whenever there exists a
constant $\lambda>0$ such that the mapping
$\lambda^{-1}\circ T$ belongs to the unit ball of the
Lipschitz space $\Lip(\Mcc,\TOM)$. In other words, the
quantity $\|\cdot\|_{LO(\Mcc,\TOM)}$ is the standard
Luxemburg norm with respect to the unit ball of
$\Lip(\Mcc,\TOM)$ (and the ``multiplication" operation
$\circ$) and $LO(\Mcc,\TOM)$ is the corresponding Orlicz
space determined by this norm, see, e.g \cite{M}. \par We
call the space $LO(\Mcc,\TOM)$ the Lipschitz-Orlicz space.
We use it to define a second ``norm": given a mapping
$T(z)=(P_z,Q_z),~z\in \Mcc,$ we put
$$
\|T\|^*:= \sup\{|D^{\alpha+\beta}
P_z(x_{Q_z})|\,r_{Q_z}^{|\beta|}:~z\in \Mcc,\, r_{Q_z}\le
1,\, |\alpha|\le k,\,|\beta|\le L-|\alpha|\},
$$
and
\bel{TS} \|T\|_{\BLO(\Mcc,\TOM)}:=
\|T\|^*+\|T\|_{LO(\Mcc,\TOM)}. \ee
We let  $\BLO(\Mcc,\TOM)$ denote the subspace of
$LO(\Mcc,\TOM)$ of ``bounded" Lipschitz mappings defined by
the finiteness of the ``norm" \rf{TS}.
\par Given closed subset $S\subset\RN$ we consider the
family of cubes $\KS$ (i.e, the cubes centered at $S$) as a
metric space equipped with the metric $\rom$, i.e., as a
subspace of the metric space $\KRN_\omega=(\KRN,\rom)$. Now
Theorem \reff{T14.2N} and Theorem \reff{T14.7N} imply the
following
\begin{theorem}\lbl{TISO}  (a). For every function
$F\in\CKLM$ there exists a mapping
$$
T(Q)=(P_Q,Q),~~~Q\in\KS,
$$
from $\BLO(\KS,\TOM)$ with $\|T\|_{\BLO(\KS,\TOM)}\le
C\|F\|_{\CKLM}$ such that $T^k_x(P_Q)=T^k_x(F)$ for every
cube $Q=Q(x,r)\in \KS$.
\par (b). Conversely, let $\omega\in\Omega_m$ be a
quasipower function. Assume that a mapping
$T(Q)=(P_Q,Q),~Q\in\KS,$ belongs to $\BLO(\KS,\TOM)$. Then
there exists a function $F\in\CKLM$ with $\|F\|_{\CKLM}\le
C\|T\|_{\BLO(\KS,\TOM)}$ such that for all $x\in S$
\bel{L-TS}
T^k_{x}(F)=\lim_{x_Q=x,\,r_Q\to 0}T^k_x(P_Q).
\ee
Moreover, for every $Q=Q(x,r)\in\KS$ and $\alpha$,
$|\alpha|\le k$,
\bel{DTF1} |D^{\alpha}T^k_{x}(F)(x)-D^{\alpha}P_{Q}(x)| \le
C\|T\|_{\BLO(\KS,\TOM)}\,r^{k-|\alpha|}\omega(r). \ee
\par Here  $C$ is a constant depending only on $k,m,n$ and
the constant $C_\omega$.
\end{theorem}
\SECT{5. Lipschitz selections of polynomial-set valued
mappings}{5}
\indent \par The ideas and results presented in Section 4
show that even though Whitney's problem for $\CKLM$ deals
with restrictions of $k$-times differentiable functions, it
is also a problem about Lipschitz mappings defined on
subsets of $\KRN$ and taking values in a very non-linear
metric space $\TOM=(\PKR,\dm)$. More specifically, the
Whitney problem can be reformulated as a problem about {\it
Lipschitz selections of set-valued mappings} from $\KS$
into $2^{\TOM}$.
\par We recall some relevant definitions: Let
$X=(\Mcc,\rho)$ and $Y=(\Tc,d)$ be metric spaces and let
$\Gc:\Mcc \to 2^{\Tc}$ be a set-valued mapping, i.e., a
mapping which assigns a {\it subset} $\Gc(x)\subset\Tc$ to
each $x\in\Mcc$. A function $g:\Mcc\to\Tc$ is said to be a
{\it selection} of $\Gc$ if $g(x)\in \Gc(x)$ for all
$x\in\Mcc$. If a selection $g$ is an element of $\Lip(X,Y)$
then it is said to be a {\it Lipschitz selection} of the
mapping $\Gc$.  (For various results and techniques related
to the problem of the existence of Lipschitz selections in
the case where $Y=(\Tc,d)$ is a Banach space, we refer the
reader to \cite{S2,S3,S4} and references therein.)
\par In \cite{F4} C. Fefferman considered the following
version of the Whitney problem: Let  $\{G(x):x\in S\}$ be a
family of convex centrally-symmetric subsets of $\PK$.
\par {\it How can we decide whether there exist $F\in\CKO$
and a constant $A>0$ such that $T_{x}^{k}(F)\in
A\circledcirc G(x)$ for all $x\in S?$}  Here $A\circledcirc
G(x)$ denotes the dilation of $G(x)$ with respect to its
center by a factor of $A$.
\par Let $P_x\in\PK$ be the center of the set $G(x)$. This
means that $G(x)$ can be represented in the form
$G(x)=P_x+\sigma(x)$ where $\sigma(x)\subset\PK$ is a
convex family of polynomials which is centrally symmetric
with respect to $0$. It is shown in \cite{F4} that, under
certain conditions on the sets $\sigma(x)$ (the so-called
condition of Whitney's $\omega$-convexity), the finiteness
property holds. The approach described in Section 4, see
Example \reff{EX-COM}, and certain ideas related to
Lipschitz selections in Banach spaces \cite{S3}, allows us
to give an upper bound for a finiteness number in
Fefferman's theorem \cite{F4}:
$$
N(k,n)=\,2^{\min\{\ell+1,\,\dim\PK\}},
$$
where $\ell=\max_{x\in S}{\dim \sigma(x)}$, see \cite{S5}.
\par This improvement of the finiteness number follows from
Fefferman's result \cite{F4} and the following
\begin{theorem}\lbl{WEAK-F}(\cite{S5}) Let $G$ be a mapping
defined on a finite set $S\subset\RN$ which assigns a
convex set of polynomials $G(x)\subset\PK$ of dimension at
most $\ell$ to every point $x$ of $S$. Suppose that, for
every subset $S'$ of $S$ consisting of at most
$2^{\,\min\{\ell+1,\,\dim\PK\}}$ points, there exists a
function $F_{S'}\in \CKO$ such that $\|F_{S'}\|_{\CKO}\leq
1$ and $T_{x}^{k}(F_{S'})\in G(x)$ for all $x\in S'$. Then
there is a function $F\in\CKO$, satisfying
$\|F\|_{\CKO}\le\gamma$ and
$$ T_{x}^{k}(F)\in G(x) ~~~for~all~~~x\in S. $$
Here $\gamma$ depends only on $k,n$ and $\card S$.
\end{theorem}
\par We use the rather informal and imprecise terminology
``$\CKO$ has the weak finiteness property" to express the
kind of result where $\gamma$ depends on the number of
points of $S$. The weak finiteness property also provides
an upper bound for the finiteness constant whenever the
strong finiteness property holds. For instance, Fefferman's
theorems in \cite{F4} reduce the problem to a set of
cardinality at most $N(k,n)$ while the weak finiteness
property decreases this number to
$2^{\min\{l+1,\dim\PK\}}$.
\par In \cite{S5}, Theorem 1.10, we show that, in turn, the
``weak finiteness" theorem, is equivalent to a certain
Helly-type criterion for the existence of a certain
Lipschitz selection of the set-valued mapping
$\Gc(x)=(G(x),x),x\in S$. An analog of this result for
set-valued mappings from $\Kc_\omega(S):=(\KS,\rom)$ into
$2^{\PKR}$ where $m=1$ and $\omega(t)=t$, is presented in
Theorem \reff{HC-OTM} below.
\par Let us see how these ideas and results can be
generalized for the space $\CKLM$ with $m>1$, and what kind
of difficulties appear in this way. We will consider the
following general version of the problem raised in
\cite{F4}.
\begin{problem}\lbl{LS-CKLM} Let  $\{G(x):x\in S\}$ be a
family of convex closed subsets of $\PK$. How can we decide
whether there exist a function $F\in\CKLM$ such that
$$ T_{x}^{k}(F)\in  G(x)
 ~~~~{\rm {\it for~all}}~~x\in S~?
$$
\end{problem}
\par In particular, if $G(x)=\{P\in \PK:~P(x)=f(x)\}, x\in
S,$ where $f$ is a function defined on $S$, this problem is
equivalent to the Whitney Problem \reff{WP-CKLM} for
$\CKLM$.
\par First, let us show that Problem \reff{LS-CKLM} is
equivalent to a Lipschitz selection problem for set-valued
mappings from $\KS_\omega:=(\KS,\rom)$ into a certain
family of subsets of $\TOM:=(\PKR,\dm)$. To this end, given
a cube $Q=Q(x,r)\in\KS$ and $\lambda>0$ we let
$H_\lambda(Q)$ denote the set of all polynomials $P\in\PL$
satisfying the following condition:
\par {\it There exists $\tP\in G(x)$ such that for every
$\alpha$, $|\alpha|\le k$,}
\bel{DAH} |D^{\alpha}\tP(x)-D^{\alpha}P(x)| \le
\lambda\,r^{k-|\alpha|}\omega(r). \ee
\par In particular,
\bel{HZERO}
H_0(Q):=\{P\in\PL:~T^k_x(P)\in G(x)\},
~~~Q=Q(x,r)\in\KS.
\ee
Clearly, $H_\lambda(Q)$ is a convex closed subset of $\PL$.
\par By $\Hc_\lambda$ we denote the set-valued mapping from
$\KS$ into $2^{\PKR}$ defined by the following formula:
$$
\Hc_\lambda(Q):=(H_\lambda(Q),Q),~~~Q\in\KS.
$$
\begin{theorem}\lbl{LS-WPL} (a) Suppose that $F\in\CKLM$
and $T^k_x(F)\in G(x)$ for every $x\in S$. Then the
set-valued mapping $\Hc_0$ has a selection
$T\in\BLO(\KS,\TOM)$ with $\|T\|_{\BLO(\KS,\TOM)}\le
C\|F\|_{\CKLM}$.
\par (b) Conversely, let $\omega\in\Omega_m$ be a
quasipower function. Suppose that there exists $\lambda>0$
such that $\Hc_\lambda$ has a selection
$T\in\BLO(\KS,\TOM)$ with
$\|T\|_{\BLO(\KS,\TOM)}\le\lambda$. Then there exists
$F\in\CKLM$ with $\|F\|_{\CKLM}\le C\lambda$ such that
$T^k_x(F)\in G(x)$ for every $x\in S$.
\par Here  $C$ is a constant depending only on $k,m,n$ and
the constant $C_\omega$.
\end{theorem}
\par {\it Proof.} (a) By part (a) of Theorem \reff{TISO}
there exists a mapping $T:\KS\to \PKR$ of the form
$T(Q)=(P_Q,Q),~Q\in\KS,$ satisfying the following
conditions: $T\in\BLO(\KS,\TOM)$,
$\|T\|_{\BLO(\KS,\TOM)}\le C\|F\|_{\CKLM}$, and
$T^k_x(P_Q)=T^k_x(F)$ for every cube $Q(x,r)\in\KS$.
\par But $T^k_x(F)\in G(x)$ so that $T^k_x(P_Q)\in G(x)$ as
well, proving that $P_Q\in H_0$, see \rf{HZERO}. Thus, the
mapping $T$ is a selection of the set-valued mapping
$\Hc_0$, and part(a) of the theorem is proved.
\par (b) Since $T$ is a selection of $\Hc_\lambda$, it can
be written in the form
$$
T(Q)=(P_Q,Q),~Q\in\KS,~P_Q\in\PL.
$$
By part Theorem \reff{TISO}, part (b), there exists
$F\in\CKLM$ with $\|F\|_{\CKLM}\le C\lambda$ such that
\rf{L-TS} and \rf{DTF1} are satisfied.
\par Prove that $T^k_x(F)\in G(x), x\in S.$ Since $T$ is a
selection of $\Hc_\lambda$, for each cube $Q=Q(x,r)\in\KS$
we have $P_Q\in H_\lambda(Q)$, so that, by \rf{DAH}, there
exists $\tP\in G(x)$ such that
$$
|D^{\alpha}\tP(x)-D^{\alpha}P_Q(x)| \le
\lambda\,r^{k-|\alpha|}\omega(r),~~|\alpha|\le k.
$$
We have
\be |D^{\alpha}T^k_x(F)(x)-D^{\alpha}\tP(x)| &\le&
|D^{\alpha}T^k_x(F)(x)-D^{\alpha}P_Q(x)|
+|D^{\alpha}P_Q(x)-D^{\alpha}\tP(x)|\nn\\
&\le&
|D^{\alpha}T^k_x(F)(x)-D^{\alpha}P_Q(x)|+
\lambda\,r^{k-|\alpha|}\omega(r),~~|\alpha|\le
k.
\nn\ee
Since $\|T\|_{\BLO(\KS,\TOM)}\le\lambda$, by \rf{DTF1},
\bel{LTIN} |D^{\alpha}T^k_x(F)(x)-D^{\alpha}\tP(x)| \le
C\lambda\,r^{k-|\alpha|}\omega(r)~~~{\rm for~ all}
~~|\alpha|\le k. \ee
\par Clearly, given $x\in S$ the quantity
$\|P\|:=\max\{|D^\alpha P(x):~|\alpha|\le k\}$ presents an
equivalent norm on the finite dimensional space $\PK$. By
inequality \rf{LTIN}, the distance from $T^k_x(F)$ to
$G(x)$ in this norm tends to $0$ as $r\to 0$. Since $G(x)$
is closed, $T^k_x(F)\in G(x)$ proving the theorem.\bx
\par As we have noted above, a geometrical background of
the weak finiteness property, Theorem \reff{WEAK-F}, is a
certain Helly-type criterion for the existence of a
Lipschitz selection. Let us formulate a version of this
result for set-valued mappings defined on finite families
of cubes in $\RN$.
\begin{theorem}\lbl{HC-OTM} Let $m=1$ and let
$\omega\in\Omega_1$ be a quasipower function. Let
$K\subset\Kc$ be a finite set of cubes in $\RN$ and let
$\Hc(Q)=(H(Q),Q),Q\in K,$ be a set-valued mapping such that
for each $Q\in K$ the set $H(Q)\subset\PK$ is a convex set
of polynomials of dimension at most $\ell$. Suppose that
there exists a constant $A>0$ such that, for every subset
$K'\subset K$ consisting of at most
$2^{\min\{\ell+1,\dim\PK\}}$ elements, the restriction
$\Hc|_{K'}$ has a selection $h_{K'}\in LO(K',\TOM)$ with
$\|h_{K'}\|_{LO(K',\TOM)}\le A$.
\par Then $\Hc$, considered as a map on all of $K$, has a
Lipschitz selection $h\in LO(K,\TOM)$ with
$\|h\|_{LO(K,\TOM)}\le \gamma A$. Here the constant
$\gamma$ depends only on $k,n$ and $\card K$.
\end{theorem}
\par We recall that for $m=1$ we have $L:=k+m-1=k$ so that
$\TOM:=(\PKR,\dm)=(\PK,\dm)$. The proof of this result
follows precisely the same scheme as in \cite{S5}.
\par It would be very useful to have such a criterion for
{\it arbitrary} $m>1$ which, in view of Theorem
\reff{LS-WPL}, would immediately lead to the weak
finiteness property for the space $\CKLM$. However, the
straightforward application of the method of proof given in
\cite{S5} to this case meets certain difficulties. In
particular, one of the crucial ingredients of the proof in
\cite{S5} is  ``consistency" of the metrics
$\rho_1(x,y):=\|x-y\|$ and  $\rho_2(x,y):=\omega(\|x-y\|)$
in the following sense: for every $x_1,x_2,x_3,x_4\in\RN$
the inequality $\rho_1(x_1,x_2)\le\rho_1(x_3,x_4)$ imply
the inequality $\rho_2(x_1,x_2)\le\rho_2(x_3,x_4)$ (and
vise versa).
\par For instance, a corresponding analog of this property
for the space $\ZM$, see \rf{NR-ZM}, is ``consistency" of
the distance
$$
\rho_1(Q_1,Q_2):=\max\{r_1,r_2\}+\|x_1-x_2\|,
~~Q_i=Q(x_i,r_i),
i=1,2,
$$
defined on the family $\Kc$ of all cubes in $\RN$, and the
metric $\rho$ defined by \rf{MHP}. However, in general,
such a ``consistency" does not hold. For example, consider
the family $\{Q_i=Q(0,r_i), i=1,2,...\}$ of cubes in $\RN$
with $r_i=2^{-i^2}$. Clearly,
$\rho_1(Q_i,Q_{i+1})=2^{-i^2}\to 0$  while
$\rho_2(Q_i,Q_{i+1})=\ln(1+r_i/r_{i+1}) \to +\infty$ as
$i\to \infty$.
\par This simple example shows that the set
$S=\{0\}\cup\{x_i,~i=1,2,...\}$ where $x_i$ are points in
$\RN$ with $\|x_i\|=2^{-i^2}$, could play a role of a
test-set in proving the weak finiteness property for the
space $\ZM$.
\renewcommand {\refname} {\centerline{\normalsize
{\bf References}}}

\end{document}